\documentclass[12pt]{amsart}
\usepackage{amsmath,amssymb,txfonts}
\usepackage{times}
\usepackage{bbm}
\usepackage[dvips]{graphicx}
\usepackage{color}
\usepackage{geometry}
\usepackage{texdraw}
\usepackage{tikz}
\allowdisplaybreaks

\usepackage[cp1252]{inputenc}

\usepackage[active]{srcltx}

\numberwithin{equation}{section}

\newtheorem{theorem}{Theorem}[section]
\newtheorem{remark}[theorem]{Remark}
\newtheorem{corollary}[theorem]{Corollary}
\newtheorem{proposition}[theorem]{Proposition}

\newtheorem{lemma}[theorem]{Lemma}

\newcommand{\ct}[1]{\langle {#1}\rangle \lower.3ex\mathrm{$_{t}$}}
\newcommand{\lt}[1]{[ {#1}] \lower.3ex\mathrm{$_{t}$}}

\def\dint{\displaystyle\int}
\def\dsum{\displaystyle\sum}
\def\dsup{\displaystyle\sup}

\numberwithin{equation}{section}

\begin{document}

	\title[Quantitative weighted bounds for the $q$-variation]
	{Quantitative weighted bounds for the $q$-variation of\\
singular integrals with rough kernels}
	
	\author{Yanping Chen}
	\address{School of Mathematics and Physics, University of Science and Technology Beijing, Beijing 100083, China}
	\email{yanpingch@ustb.edu.cn}

\author{Guixiang Hong}
\address{School of Mathematics and Statistics, Wuhan University,
	Wuhan 430072, China}
	\email{guixiang.hong@whu.edu.cn}
	
	\author{Ji Li}
\address{Department of Mathematics,
		Macquarie University, NSW, 2109, Australia}
	\email{ji.li@mq.edu.au}

	\thanks{The project was in part supported by: Yanping Chen's
		NSF of China (\# 11871096, \# 11471033); Guixiang Hong's NSF of China ( \# 11601396); Ji Li's  ARC DP 170101060.}

	
	\date{\today}
	

	\begin{abstract}
		In this paper, we study the quantitative weighted bounds for the $q$-variational singular integral operators with rough kernels. The main result is
for the sharp truncated singular integrals itself
		$$		\|V_q\{T_{\Omega,\varepsilon}\}_{\varepsilon>0}\|_{L^p(w)\rightarrow L^p(w)}\leq c_{p,q,n} \|\Omega\|_{ L^\infty}(w)_{A_p}^{1+1/q}\{w\}_{A_p},$$
where the quantity $(w)_{A_p}$, $\{w\}_{A_p}$ will be recalled in the introduction;
		 we do not know whether this is sharp, but it is the best known quantitative result for this class of operators, since when $q=\infty$, it
coincides with the best known quantitative bounds by Di Pilino--Hyt\"{o}nen--Li or Lerner. In the course of establishing the above estimate, we obtain several
quantitative weighted bounds which are of independent interest. We hereby highlight two of them. The first one is
		 $$		\|V_q\{\phi_k\ast T_{\Omega}\}_{k\in\mathbb Z}\|_{L^p(w)\rightarrow L^p(w)}\leq c_{p,q,n}  \|\Omega\|_{ L^\infty}(w)_{A_p}^{1+1/q}\{w\}_{A_p},$$
		 where $\phi_k(x)=\frac1{2^{kn}}\phi(\frac x{2^k})$ with $\phi\in C^\infty_c(\mathbb R^n)$ being any non-negative radial function, and the sharpness for $q=\infty$ is due to
Lerner; the second one is
		  $$		\|\mathcal{S}_q\{T_{\Omega,\varepsilon}\}_{\varepsilon>0}\|_{L^p(w)\rightarrow L^p(w)}\leq c_{p,q,n}  \|\Omega\|_{ L^\infty}(w)_{A_p}^{1/q}\{w\}_{A_p},$$
		  and the sharpness for $q=\infty$ follows from the Hardy--Littlewood maximal function.

	\end{abstract}
	
	\maketitle
	
\section{Introduction and statement of main results}

\subsection{Background and main result}

Let $\mathcal F=\{F_t:t\in\mathcal I\subset \mathbb{R}\}$ be a family of Lebesgue measurable functions defined on $\mathbb{R}^n$. For  $x$ in
$\mathbb{R}^n$, the value of the  $q$-variation function $V_q(\mathcal F)$ of the family $\mathcal F$ at $x$ is defined by
\begin{align}\label{VqF}
V_q(\mathcal F)(x):= \sup\Big(\sum_{k\geq1}
|F_{t_k}(x)-F_{t_{k-1}}(x)|^q\Big)^{\frac{1}{q}} ,\quad q\ge1,
\end{align}
where the supremum runs over all increasing subsequences $\{t_k:k\geq0\}\subset \mathcal I$.
Suppose $\mathcal{A}=\{{A}_t\}_{t\in\mathcal I}$ is a family of operators on $L^p(\Bbb R^n)\, (1\le p\le\infty)$, the associated strong $q$-variation operator
$V_q\mathcal A$ is defined as
$$V_q\mathcal A(f)(x)=V_q(\{{A}_tf(x)\}_{t\in\mathcal I}).$$

\smallskip
There are two elementary but important observations that motivate the development of variational inequalities in ergodic theory and harmonic analysis. The
first one is that from the fact that $V_q(\mathcal F)(x)<\infty$ with finite $q$ implies the convergence of $F_t(x)$ as $t\rightarrow t_0$ whenever $t_0$
is an adherent point for $\mathcal I$, it is easy to observe that $A_t(f)$ converges a.e.  for $f\in L^p$ whenever $V_q\mathcal A$ for finite
$q$ is of weak type $(p,p)$  with $p<\infty$.  The second one is that $q$-variation function dominate pointwisely the maximal function: for any $q\ge1$,
$$
A^\ast(f)(x)\le A_{t_0}f(x)+ V_{q}(\mathcal{A}f)(x),
$$
where $A^\ast$ is the maximal operator defined by
$A^\ast(f)(x):=\sup_{t\in\mathcal I}|{A}_t(f)(x)|$ and $t_0\in\mathcal I$ is any fixed number.

\smallskip

Because of the first observation, Bourgain \cite{Bou89} introduced in ergodic theory the first variational estimate which had originated from the
regularity of Brownian motion in martingale theory \cite{Lep76, PiXu88} since the pointwise convergence may not hold for any non-trivial dense subclass of
functions in some ergodic models and thus the maximal Banach principle does not work, see \cite{JKRW98, JRW03} for further results in this direction. Because of
the second observation, many maximal inequalities in harmonic analysis, such as maximal singular integrals, maximal operators of Radon type, the
Carleson--Hunt theorem, Cauchy type integrals, dimension free Hardy--Littlewood maximal estimates etc,  have been strengthened to variational inequalities
\cite{CJRW2000, CJRW2002, JSW08, MiTr14,  MST15, OSTTW12, MaTo, BMSW, Zor14, MTZ14}.  There are also many other publications on vector-valued and weighted norm estimates coming to
enrich the literature on this subject (cf. e.g. \cite{GiTo04,LeXu2,DMT12,JSW08,Hon1}).
\smallskip
In this paper, we continue to study $q$-variational estimates but focus on the quantitative weighted bounds for singular integrals with rough kernels.

\smallskip
For $\varepsilon>0,$ suppose that $T_{\Omega,\varepsilon}$ is the truncated homogeneous singular integral operator defined by
\begin{align}\label{ts}
T_{\Omega, \varepsilon}
f(x)=\int_{|y|>\varepsilon}\frac{\Omega(y')}{|y|^n}f(x-y)dy,
\end{align}
where $y'={y/|y|}$, $\Omega$ is
homogeneous of degree zero, integrable on $\mathbb S^{n-1}$  and
 satisfies the cancellation condition
\begin{align}\label{zero}
\int_{\mathbb S^{n-1}}\Omega(y')d\sigma(y')=0.
\end{align}
 Denote the family $\{T_{\Omega, \varepsilon}\}_{\varepsilon>0}$ of operator by $\mathcal T_\Omega$. For $1\leq p<\infty$ and $f\in L^p(\mathbb R^n)$,
 the Calder\'on--Zygmund singular integral operator $T_\Omega$ with homogeneous kernel could be defined by
\begin{align}\label{dense}
T_\Omega f(x)=\lim_{\varepsilon\rightarrow0^+}T_{\Omega, \varepsilon} f(x), \ a.e. \ x\in\mathbb R^n.
\end{align}
 This definition is justified by the boundedness of the maximal singular integral $T_\Omega^\ast$ which is given by $T_\Omega^\ast
 f(x)=\sup_{\varepsilon>0}|T_{\Omega, \varepsilon} f(x)|$ together with the maximal Banach principle.

  \smallskip
  In 2002, Campbell {\it et al} \cite{CJRW2002} first proved the
 variational inequalities for $  T_\Omega $ extending their previous results for the Hilbert transform \cite{CJRW2000},  which gives extra information on
 the speed of convergence in \eqref{dense}. See also \cite{JSW08, DHL} for more results.

\smallskip\smallskip\par
\noindent {\bf Theorem A.}\quad {\it Let $\Omega\in L\log^+\!\!L(\mathbb S^{n-1})$ satisfying \eqref{zero} and $\mathcal T_\Omega$ be as defined above. Then the variational inequality
\begin{align*}
\|V_q\mathcal T_\Omega (f)\|_{L^p}\le c_{p,q,n}\|\Omega\|_{L\log^+\!\!L(\mathbb S^{n-1})}\|f\|_{L^p}
\end{align*}
holds for $1<p<\infty$ and $q>2$.}
\par

\smallskip

The first weighted norm $q$-variational inequality, to the authors' best knowledge,  is due to Ma, Torrea and Xu in \cite{MTX1}. Later on, their results were extended to higher dimensions \cite{MTX2, KZ14}. The weighted version of Theorem A was also considered in \cite{CDHL}. On the other hand, motivated by the solution of the $A_2$ conjecture,
 Hyt\"{o}nen,  Lacey and P\'{e}rez \cite{HLP} established the sharp weighted inequality for the variation of
 the smooth truncations of Calder\'{o}n--Zygmund operator with $\omega$-Dini type regularity;  in \cite{deZ16}, they  extended
the sharp weighted estimates for smooth truncations
to the case of sharp truncations. The latter result was generalized to matrix weight by Duong, Li and Yang \cite{DLY19}.
\smallskip
However, the (sharp) quantitative weighted version of Theorem A, that is for operators without regularity,  has not been considered up to now except the
special case $q=\infty$ due to Lerner \cite{Ler 2} and  Di Pilino, Hyt\"{o}nen and Li \cite{PHL}.
\smallskip

One main goal of the present paper is to address this question, and to provide the quantitative weighted variational inequality for $\mathcal T_\Omega$.
\quad

\begin{theorem}  \label{thm5}  Let $\Omega\in L^\infty({\Bbb S}^{n-1})$ satisfying \eqref{zero} and $\mathcal T_\Omega$ be as defined above. Then the variational inequality
\begin{align}\label{l1 thm5}
\|V_q\mathcal T_\Omega (f)\|_{L^p(w)}\le c_{p,q,n}\|\Omega\|_{ L^\infty}(w)_{A_p}^{1+1/q}\{w\}_{A_p}\|f\|_{L^p(w)}
\end{align}
holds for $q>2$, $1<p<\infty$ and $w\in A_p$, where the implicit constant $c_{p,q,n}$ is independent of $f$ and $w$, and
$(w)_{A_p}:=\max\Big\{[w]_{A_\infty},\,[w^{1-p'}]_{A_\infty}\Big\}$, $\{w\}_{A_p}:=[w]_{A_p}^{1\over p}\max\Big\{[w]_{A_\infty}^{1\over
p'},\,[w^{1-p'}]_{A_\infty}^{1\over p}\Big\}$, and the constant $[w]_{A_p}$, $1<p\leq\infty$, is the norm of $w$, given in \eqref{Ap} and \eqref{Ainfty}.
\end{theorem}

\begin{remark}
Recall that the
(sharp) quantitative weighted bounds  for $ T_\Omega$ $($given as
in  \eqref{ts} with $\Omega\in L^\infty({\Bbb S}^{n-1})$ satisfying \eqref{zero}$) $ have been studied intensively in the last three years $($see for example
\cite{HRT,CCDO,Ler 0, Ler 2}$)$ with the key tool---sparse domination whose pointwise version was originated in \cite{Lac15}. Our quantitative estimate in
\eqref{l1 thm5} for $q=\infty$ matches the best known results in \cite{Ler 2, PHL} for $T_\Omega^*$. However we have no idea whether the quantity on the
right hand side of \eqref{l1 thm5} could be independent of $q$, let alone the sharpness of the weighted constant.

\end{remark}
\quad
\subsection{Approach and main methods} We now state the methods and techniques for proving our main result. Some estimates are of independent interest, and
hence we will single them out as theorems or propositions.

\smallskip


As usual, we shall prove 
Theorem \ref{thm5}  by verifying separately the corresponding inequalities for the long and short variations as follows  (see e.g. \cite[Lemma 1.3]{JSW08}). 

\quad

\begin{theorem}\label{pro01} Let $q>2.$
Let $\Omega\in L^\infty({\Bbb S}^{n-1})$ satisfying \eqref{zero}.  Then for $1<p<\infty$ and $w\in A_p,$
\begin{align}\label{long} \|V_q(\{T_{\Omega,2^k}f\}_{k\in \Bbb Z})\|_{L^p(w)}\le c_{p,q,n}\|\Omega\|_{L^\infty}(w)_{A_p}^{1+1/q}\{w\}_{A_p}\|f\|_{L^p(w)}.
\end{align}
\end{theorem}
\quad

\begin{theorem}\label{pro02}
Let $q\geq2.$
Let $\Omega\in L^\infty({\Bbb S}^{n-1})$ satisfying \eqref{zero}. Then for $1<p<\infty$ and $w\in A_p,$
\begin{align}\label{short} \|\mathcal{S}_q({\mathcal{T}}_\Omega f)\|_{L^p(w)}\le c_{p,q,n}\|\Omega\|_{L^\infty}(w)_{A_p}^{1/q}\{w\}_{A_p}\|f\|_{L^p(w)}, \end{align}
where
$$\mathcal{S}_q({\mathcal{T}}_\Omega f):=\bigg(\sum_{k\in\mathbb Z}[V_{q}(\{{ T}_{\Omega,\varepsilon} f\}_{2^k\leq\varepsilon<2^{k+1}})]^q\bigg)^{1/q}.$$
\end{theorem}

\quad

The proof of Theorem \ref{pro01} is given in Section 6 and the proof of Theorem \ref{pro02} is given in Section 7, 
while the important tools for the proof of these two theorems are established in Sections 3, 4 and 5.

 It is worthy to point out the above two results are of independent interest due to the presence of different quantitative behavior.  In particular, the
 estimate in \eqref{short} is sharp in the case $q=\infty$ which follows from the sharp weighted norm estimate of the Hardy-Littlewood maximal operator.

To help understand these
structures in the methods and techniques, we can look at the following flow chart where all the terms are listed.

\begin{center}
\begin{tikzpicture}
\node[draw,rounded corners] (1) at (0,16.2){${\rm Thm} 1.1:\ \ V_q\mathcal T_\Omega (f)$};
\node[draw,rounded corners] (2) at (6,15){${\rm Thm} 1.3$};
\node[draw,rounded corners] (3) at (-5.5,15){${\rm Thm} 1.4$};
\node[draw,rounded corners] (4) at (-5,13){$\mathcal{T}_\Omega^1(f):\ \phi_{k}\ast T_\Omega$};
\node[draw,rounded corners] (5) at (0.4,13){$\mathcal{T}_\Omega^2(f):\ \phi_{k}\ast (K_\Omega\chi_{|\cdot|\leq 2^k})\ast f$};
\node[draw,rounded corners] (6) at (6,13){$\mathcal{T}_\Omega^3(f):\ (\delta_0-\phi_{k})\ast T_{\Omega,2^k}f$};
\node[draw,rounded corners] (7) at (-6,8){${\rm Thm} 1.6$};
\node[draw,rounded corners] (8) at (-3,11){${\rm Thm} 1.5:\ V_q\Phi( T_\Omega f)$};
\node[draw,rounded corners] (10) at (2,10){${\rm Prop} 1.8:\ \phi_k*K_k$};
\node[draw,rounded corners] (11) at (5,7){${\rm Prop} 1.9: \ ({{\phi}}_k-\delta_0)\ast K^k$};
\draw (1) to (2);
\node at (-3,15.5) {Short Variation};
\draw (1) to (3);
\node at (3,15.5) {Long Variation};
\draw (2) to (5);
\node at (3.7,14) {Cotlar\ \ \ type\ \ \  decomposition};
\draw (2) to (4);
\draw (2) to (6);
\draw (4) to (8);
\draw (5) to (10);
\draw (6) to (11);
\draw (7) to (8);
\node at (-4,9.5){Dini to Rough kernel};
\draw (7) to (10);
\draw (7) to (11);
\end{tikzpicture}
\end{center}

\medskip\noindent
\quad
\subsubsection{{\bf Methods for proving Theorem \ref{pro01}}:}As a standard way to study variational inequalities for rough singular integrals  (see e.g.
\cite{CDHL,DHL} originating from \cite{DR86}), the first step is to exploit the Cotlar type decomposition of $T_{\Omega,2^k}$. Let $\phi\in
C^{\infty}_{0}(\mathbb{R}^n)$ be a non-negative radial function, supported in $\{|x|\leq1/4\}$ with $\int\phi(x)dx=1$, and denote
$\phi_{k}(x):= 2^{-kn}\phi(2^{-k}x)$. {Let $K_\Omega(x):={\Omega(x/|x|)}{|x|^{-n}}$, and $\delta_0$ be the
Dirac measure at $0$.
Then for a reasonable function
$f$,
we decompose $T_{\Omega,2^k}f$ as
\begin{align}\label{rough dec}
T_{\Omega, 2^k}f
&=\phi_{k}\ast T_\Omega f-\phi_{k}\ast (K_\Omega\chi_{|\cdot|\leq 2^k})\ast f+(\delta_0-\phi_{k})\ast T_{\Omega,2^k}f,
\end{align}} which divides $\mathcal{T}_\Omega =\{T_{\Omega, 2^k}\}_{k\in \Bbb Z}$ into three families:
\begin{align*}
\mathcal{T}_\Omega^1(f)&:=\{\phi_{k}\ast T_\Omega f\}_{k\in \Bbb Z},\\
\mathcal{T}_\Omega^2(f)&:=\{\phi_{k}\ast (K_\Omega\chi_{|\cdot|\leq 2^k})\ast f\}_{k\in \Bbb Z},\\
\mathcal{T}_\Omega^3(f)&:=\{(\delta_0-\phi_{k})\ast T_{\Omega,2^k}f\}_{k\in \Bbb Z},
\end{align*}
Thus, it suffices to verify  the weighted $L^p$ estimate for the variation of $\mathcal{T}_\Omega^i(f)$,
$i=1,2,3$.
%

The second usual step to deal with variational estimates for rough singular integrals is to exploit the almost orthogonality principle based on
Littlewood-Paley decomposition and interpolation. However, this step is not obvious at all in the present setting for quantitative weighted estimates.
Indeed, the standard Littlewood-Paley decomposition seems to be insufficient. We will exploit the one associated to a sequence of natural numbers
$\mathcal N=\{N(j)\}_{j\geq0}$ used in \cite{HRT}; moreover, it is quite involved to gain the sharp kernel estimates involving the smoothing version of
$T_\Omega$ which are necessary to obtain the sharp quantitative weighted estimates. Now let us comment on the proof term by term in details.

\medskip

 {\bf Comments on the estimate of $\mathcal{T}_\Omega^1(f)$.} Let us first formulate below the desired estimate of this term which is of independent interest since this
 variational estimate strengthens the result of Lerner \cite{Ler 2}, and is sharp in the case $q=\infty$ according to \cite{Ler 2}.
Denote by $\Phi(g)=\{\phi_j\ast g\}_{j\in \Bbb Z}.$


\begin{theorem}  \label{thm3}  Let $q>2$. Let $\phi$ be as used  in \eqref{rough dec} and $ T_\Omega $ be given as
in  \eqref{ts} with $\Omega\in L^\infty(S^{n-1})$ satisfying \eqref{zero}. Then for $1<p<\infty$ and $w\in A_p,$
\begin{align}\label{l1t}
\|V_q\Phi( T_\Omega f) \|_{L^p(w)}\le c_{p,q,n} \|\Omega\|_{ L^\infty}(w)_{A_p}^{1+1/q}\{w\}_{A_p}\|f\|_{L^p(w)}.
\end{align}


\end{theorem}
The proof of Theorem \ref{thm3} is given in Section 5.


To obtain \eqref{l1t}, we apply the Littlewood-Paley decomposition of $T_\Omega$ in \cite{HRT} to get
$
T_{\Omega}=\sum_{j=0}^{\infty} {T}_{j}^{\mathcal N},
$
where each $ {T}_{j}^{\mathcal{N}}$ is a $\omega_j$-Dini Calder\'on--Zygmund operator  of convolution type.
This yields
$$
\|V_q\Phi(T_\Omega f)\|_{L^p(w)}\leq\sum_{j=0}^{\infty} \|V_q\Phi(T_{j}^{\mathcal N}f)\|_{L^p(w)}.
$$
Then to control the summation of the right-hand side of the above inequality, it suffices to show for each term the unweighed estimate with rapid decay
(with respect to $j$) and the sharp weighted norm $L^p$ estimate for all $1<p<\infty$ and then to apply the Stein--Weiss interpolation---Lemma \ref{sw}. The rapid
decay estimate will follow from the $L^2$ estimate of $T^\mathcal N_j$; the involved part lies in the sharp estimates in terms of the weighted
constant and the Dini norm. The latter is not only necessary for further application but also is of independent interest, we formulate it in
the following theorem.

\begin{theorem}  \label{thm1}
Let $q>2$. Let $T$ be a $\omega$-Dini Calder\'on--Zygmund operator of convolution type with the kernel $K$ satisfying the
following cancellation condition:
for any $0<\varepsilon,\,R<\infty,$
\begin{align}\label{2} \int_{\varepsilon<|x|<R}
K(x)\,dx=0.
\end{align}
Then for any $1<p<\infty$ and $w\in A_p,$
\begin{align}\label{phia}\|V_q\Phi(Tf)\|_{L^{p}(w)}\lesssim c_{p,q,n}(\|{\omega}\|_{Dini}+\|\omega\|_{Dini}^{1/q}\|\omega^{1/q'}\|_{Dini}+\|\omega\|_{Dini}^{1/q'}\|\omega^{1/2}\|_{Dini}^{2/q}+C_K+\|T\|_{L^2\rightarrow L^2})\{w\}_{A_p}\|f\|_{L^{p}(w)}.\end{align}
\end{theorem}
\quad

The proof of Theorem \ref{thm1} is given in Section 5.

Firstly, it would be impossible to get the sharp weighted bound in \eqref{phia} via iteration by showing the sharp weighted estimate of $V_q\Phi$ since the
one of $T$ is already linear (see. e.g. \cite{HRT, Ler 0}); let alone, it is not clear at the moment whether one can  obtain the sharp weighted bound of $V_q\Phi$.
One of the obstacles is due to
the fact that the variation norm is not monotone  in the sense that  $a_k\leq b_k\ (\forall k\in\mathbb Z)$ does not imply $V_q(\{a_k\}_k)\leq V_q(\{b_k\}_k)$ unless $q=\infty$.

We provide the proof of \eqref{phia} based on the Cotlar type decomposition
\begin{align}\label{rough dec2}
\phi_k\ast Tf
&=K^k\ast f+\phi_{k}\ast K_k\ast f+(\phi_k-\delta_{0})\ast K^k\ast f,
\end{align}
where $K_k=K\chi_{|\cdot|\le 2^{k}}$ and $K^k=K\chi_{|\cdot|\ge 2^{k}}.$
Now, the sharp weighted bound of  $V_q(\{K^k\ast f\}_{j\in \Bbb Z})$ has been obtained in \cite{deZ16} (see also \cite{DLY19, CM}) with bounds
$\|{\omega}\|_{Dini}+C_K+\|T\|_{L^2\rightarrow L^2}$. Regarding the second term, at a first glance, it should be easier to handle than the second
 term  in \eqref{rough dec} since the kernel $K$ enjoys some regularity; but it turns out that they are in the same level of
complexity as they will follow from Proposition \ref{pro2} below. Finally, the third term in \eqref{rough dec2}  is easier to deal with than  the third
 term  in \eqref{rough dec}, which will follow from Proposition \ref{pro4}. In this process, together with the weak type $(1,1)$ estimates of variational Calder\'on--Zygmund operators
\cite{DLY19}, we also obtain the following result as a byproduct.

\begin{corollary}
Let $q>2$ and $T$ be as in Theorem \ref{thm1}.
Then
\begin{align*}\|V_q\Phi(Tf)\|_{L^{1,\infty}}\lesssim c_{q,n}(\|{\omega}\|_{Dini}+\|\omega\|_{Dini}^{1/q}\|\omega^{1/q'}\|_{Dini}+\|\omega\|_{Dini}^{1/q'}\|\omega^{1/2}\|_{Dini}^{2/q}+C_K+\|T\|_{L^2\rightarrow L^2})\|f\|_{L^1}.\end{align*}
\end{corollary}

\medskip

{\bf Comments on the estimate of $\mathcal{T}_\Omega^2(f)$.} 
About $\mathcal{T}_\Omega^2(f),$ we start with the comments on the estimate of the second term in the decomposition \eqref{rough dec2} since the latter involving regular kernel $K$ seems to be easier to handle.
%
In order to apply the criterion for sharp weighted norm estimate provided by Lerner in \cite{Ler 0}, we first control the variation norm by a stronger norm. That is,
$$V_q(\{\phi_{k}\ast K_k\ast f\}_k)=\|\phi_{k}\ast K_k\ast f\|_{V_q}\leq 2\|\phi_{k}\ast K_k\ast f\|_{\ell_q}=:L_q(\{\phi_{k}\ast K_k\ast f\}_k).$$
Then all the desired estimates for the left quantity involving $V_q$ will follow from that involving $L_q$.
\smallskip
We establish the following result.

 \begin{proposition} \label{pro2}
Let $K$ be a kernel satisfying the mean value zero property \eqref{2} and the size condition \eqref{K0}. Then the following assertions hold.
 \begin{enumerate}
 \item The size condition: for $x\in \mathbb R^n\backslash\{0\}$
 \begin{align} \label{sd}
 |{\phi}_k\ast K_k(x)|\le c_n\frac{C_K}{2^{kn}},\,\,\,k\in \Bbb Z.
\end{align}
The Lipchitz condition:
 for $0<2|y|\le |x|,$
 \begin{align}\label{smd}
 \sum_{k\in \Bbb Z}|{\phi}_k\ast K_k(x-y)-\phi_k\ast K_k(x)|&\le c_nC_K\frac{|y|^\gamma}{|x|^{n+\gamma}}.\end{align}
 \item  $L^2$ boundedness:
\begin{align}\label{s2d}
\|(\sum_{k\in \Bbb Z}|{\phi}_k\ast K_k\ast f|^2)^{1/2}\|_{L^{2}}\le c_nC_K\|f\|_{L^2}.
\end{align}
Weak type $(1,1)$ boundedness:
\begin{align}\label{w1d}
\|(\sum_{k\in \Bbb Z}|{\phi}_k\ast K_k\ast f|^2)^{1/2}\|_{L^{1,\infty}}\le c_nC_K\|f\|_{L^1}.
\end{align}
Weak type $(1,1)$ boundedness of the grand maximal truncated operator   (see \eqref{m}):
\begin{align}\label{w1md}
\|\mathcal M_{(\sum_{k\in \Bbb Z}|{\phi}_k\ast K_k\ast f|^2)^{1/2}}\|_{L^{1,\infty}}\le c_nC_K\|f\|_{L^1}.
\end{align}

\item Sharp weighted norm inequalities: for all $1<p<\infty,$ and $w\in A_p,$
\begin{align}\label{swd}
\bigg\|\bigg(\sum_{k\in \Bbb Z}|\phi_k\ast K_k\ast f|^2\bigg)^{1/2}\bigg\|_{L^{p}(w)}&\le c_{p,n} C_K
\{w\}_{A_p}\|f\|_{L^{p}(w)}.
\end{align}
 \end{enumerate}

\end{proposition}


The proof of Proposition \ref{pro2} is given in Section 3.

Surely, by the fact that $\ell^2\subset\ell^q$ for $q\geq2$, the desired estimate for the second term in the decomposition \eqref{rough dec2} can be deduced from the above proposition. Moreover, the above result is somewhat surprising in the sense that there is no need of any regularity assumption on $K$; in particular, it applies well to $\mathcal{T}_\Omega^2(f)$, and thus the almost orthogonality technique can be avoided, improving the corresponding argument in \cite{CDHL}.


\medskip

{\bf Comments on the estimate of $\mathcal{T}_\Omega^3(f)$.}
Again, we start with the comments on the estimate of the third term in the decomposition \eqref{rough dec2}.  We prove the following result.
\begin{proposition}\label{pro4}
 Let $q>2.$  Let $K$ be a kernel satisfying the mean value zero property \eqref{2} and the size condition \eqref{K0} and the Dini condition \eqref{K1}.
 Then the following assertions hold.
 \begin{enumerate}
 \item The size condition: for $x\in \mathbb R^n\backslash\{0\}$
 \begin{align} \label{su}
\sum_{k\in \Bbb Z}|({{\phi}}_k-\delta_0)\ast K^k(x)|&\le c_{n} (\|\omega\|_{Dini}+C_K)\frac{1}{|x|^{n}}.
\end{align}
The Dini condition:
 for $0<2|y|\le |x|,$
 \begin{align}\label{smu}
 \bigg(\sum_{k\in \Bbb Z}|({{\phi}}_{k}-\delta_0)\ast &K^k(x-y)-({{\phi}}_{k}-\delta_0)\ast K^k(x)|^q\bigg)^{1/q}\\&\le c_{q,n}\frac{C_K
 |y|^{\theta/2}}{|x|^{n+\theta/2}}+c_{q,n}\frac{\omega(|y|/|x|)^{\frac{1}{q'}}}{|x|^{n}}\|\omega\|_{Dini}^{1/q}+c_{q,n}\sum_{{2^k}\ge
 |y|}\frac{C_K}{|x|^n}\chi_{2^{k+1}-|y|\le|x|\le2^{k+1}+|y|}.\nonumber
 \end{align}
 \item  $L^2$-boundedness:
\begin{align}\label{s2u}
\|(\sum_{k\in \Bbb Z}|(\delta_0-{\phi}_k)\ast K^k\ast f|^q)^{1/q}\|_{L^{2}}\le c_{q,n} (C_{K}+\|\omega\|_{Dini}+\|\omega\|_{Dini}^{1/q'}\|\omega^{1/2}\|_{Dini}^{2/q})\|f\|_{L^2}.
\end{align}
Weak type $(1,1)$ boundedness:
\begin{align}\label{w1u}
 \big\|(\sum_{k\in \Bbb Z}|(\phi_{k}-\delta_0)\ast K^k\ast
f|^q)^{1/q}\big\|_{L^{1,\infty}}\le  c_{q,n}(C_K+\|\omega\|_{Dini}+\|\omega\|_{Dini}^{1/q}\|\omega^{1/q'}\|_{Dini}+\|\omega\|_{Dini}^{1/q'}\|\omega^{1/2}\|_{Dini}^{2/q})\|f\|_{L^{1}}.
\end{align}
Weak type $(1,1)$ boundedness of the grand maximal truncated operator:
\begin{align}\label{w1mu}
\|\mathcal{M}_{(\sum_{k\in \Bbb Z}|(\phi_{k}-\delta_0)\ast K^k\ast f|^q)^{1/q}}\|_{L^{1,\infty}}\le
c_{q,n}(C_K+\|\omega\|_{Dini}+\|\omega\|_{Dini}^{1/q}\|\omega^{1/q'}\|_{Dini}+\|\omega\|_{Dini}^{1/q'}\|\omega^{1/2}\|_{Dini}^{2/q})\|f\|_{L^{1}}.
\end{align}

\medskip
\item Sharp weighted norm inequalities: for all $1<p<\infty,$ and $w\in A_p,$
\begin{align}\label{swu}
\big\|\big(\sum_{k\in \Bbb Z}|(\phi_k-\delta_0)\ast K^k\ast f|^q\big)^{1/q}\big\|_{L^p(w)}&\le
c_{p,q,n}(\|{\omega}\|_{Dini}+\|\omega\|_{Dini}^{1/q}\|\omega^{1/q'}\|_{Dini}+\|\omega\|_{Dini}^{1/q'}\|\omega^{1/2}\|_{Dini}^{2/q}+C_K)\{w\}_{A_p}\|f\|_{L^p(w)}.
\end{align}
 \end{enumerate}
\end{proposition}

This proposition will be shown in Section 4.

It is a surprise that at the moment of writing we are unable to show the above result for $q=2$ since our argument in showing the weak type $(1,1)$ estimate of the localized maximal operator \eqref{w1mu} depends heavily on the same estimate of $q$-variation operator associated to the Hardy--Littlewood averages.
\smallskip
Here the regularity of $K$ plays an important role in checking the size and smooth conditions \eqref{su} \eqref{smu} etc. For the rough kernel $K_\Omega$ in $\mathcal{T}_\Omega^3(f)$, we have to decompose it by using the Littlewood--Paley decomposition from \cite{HRT} and then exploit the almost orthogonality principle. 
 There are several ways to accomplish that. One natural way is
as in showing \eqref{l1t}, to decompose like $K_\Omega\chi_{|\cdot|>2^k}=\sum^\infty_{j=0}K^\mathcal N_j\chi_{|\cdot|>2^k}$, where $K^\mathcal N_j$ is the kernel of $T^\mathcal N_j$ mentioned previously.  Then for each $j$, $K^\mathcal N_j$ satisfies the assumption in the above proposition and the desired sharp weighted norm estimate follows; however, in this case, the rapid decay $L^2$ estimate is hard to gain since the kernel $K^\mathcal N_j$ is not homogeneous anymore and  Van der Corput's lemma does not apply.
\smallskip
Instead, we make use of another decomposition. Firstly write $$K_{\Omega}\chi_{|\cdot|>2^k}=\sum_{s\geq0}v_{k+s},$$
where $v_{k}(x):={\Omega(x/|x|)}{|x|^{-n}}\chi_{2^k<|x|\le 2^{k+1}}(x)$ for $k\in\mathbb Z$; then decompose each $v_{k+s}$ as
$$v_{k+s}=\sum_{j\geq0}K^\mathcal N_{k+s,j},$$
where $K^\mathcal N_{k+s,j}$ is the kernel of $T^\mathcal N_{k+s,j}$ (see \eqref{ksj} for the exact definition).
Then for each $K^\mathcal N_{k,j}=\sum_{s\geq0}K^\mathcal N_{k+s,j}$, we can show the desired unweighted norm
estimate with rapid decay in $j$ and the sharp weighted norm estimate and then the Stein--Weiss interpolation applies. See Section 6 for details.





\subsubsection{{\bf Methods for proving Theorem \ref{pro02}:}}\ \
Following from the sharp weighted boundedness of the Hardy--Littlewood maximal operator (see \eqref{s1f}),
\begin{align*}
\|\mathcal{S}_{\infty}({\mathcal{T}}_\Omega f)\|_{L^p(w)}&\le c_{p,n} \|\Omega\|_{L^\infty}\{w\}_{A_p}\|f\|_{L^p(w)}.
\end{align*}
Thus by interpolation, it suffices to prove that
\begin{align}\label{s2f in}
\|\mathcal{S}_{2}(\mathcal{T}_\Omega f)\|_{L^p(w)}&\le c_{p,n}  \|\Omega\|_{L^\infty}(w)_{A_p}^{1/2}\{w\}_{A_p}\|f\|_{L^p(w)}.
\end{align}

To verify this, we write
\begin{align}
\mathcal{S}_{2}({\mathcal{T}}_\Omega f)(x)&=\Big(\sum_{k\in\mathbb{Z}}|V_{2,k}( f)(x)|^2\Big)^{\frac{1}{2}}:
=\Big(\sum_{k\in\mathbb{Z}}\|\{T_{k,t}f(x)\}_{t\in[1,2]}\|_{V_2}^2\Big)^{\frac{1}{2}},
\end{align}
where $T_{k,t}f(x)=v_{k,t}\ast f(x)$ for $t\in [1,2]$, and $v_{0,t}(x)={\Omega(x/|x|)}{|x|^{-n}}\chi_{_{t\leq|x|\leq2}}(x)$ and
$v_{k,t}(x)={2^{-kn}}\nu_{0,t}(2^{-k}x)$ for $k\in\mathbb{Z}$.
Next, by using the Littlewood--Paley decomposition as in \cite{HRT}, we further have
$T_{k,t}=\sum_{j=0}^\infty T^\mathcal{N}_{k, t, j}$, whose definition will be given in Section 7.
Therefore, by the Minkowski inequality, we get
\begin{align*}
\mathcal{S}_2({\mathcal{T}}_\Omega f)(x)
&\leq\sum^\infty_{j=0}\mathcal{S}_{2,j}^\mathcal{N}( f)(x),
\end{align*}
where
$$
\mathcal{S}_{2,j}^\mathcal{N}( f)(x) := \bigg(\sum_{k\in\mathbb{Z}}\bigg\|\bigg\{T_{k,t,j}^\mathcal{N}f(x)\bigg\}_{t\in[1,2]}\bigg\|_{V_2}^2\bigg)^\frac{1}{2}.
$$
Hence, to prove \eqref{s2f in}, it again suffices to verify for $\mathcal{S}_{2,j}^\mathcal{N}( f)$ the unweighted norm with rapid decay in $j$ and the sharp weighted norm estimate and then to apply the Stein--Weiss interpolation---Lemma \ref{sw}.

\bigskip
\noindent{\bf Notation}. From now on, $p'=p/(p-1)$ represents the conjugate number of $p\in [1,\infty)$; $X\lesssim Y$ stands for $X\le C Y$ for a constant
$C>0$ which is independent of the essential variables living on $X\ \&\ Y$; and $X\approx Y$ denotes  $X\lesssim Y\lesssim X$.

\section{Preliminaries}

\subsection{Muckenhoupt weights}

We first recall the definition and some properties of $A_p$ weights on $\mathbb R^n$. Let $w$ be a non-negative locally integrable function defined on
$\mathbb R^n$. We say $w\in A_1$ if there is a constant $C>0$ such that $M(w)(x)\le Cw(x)$, where $M$ is the Hardy-Littlewood maximal function.
Equivalently, $w\in A_1$ if and only if there is a constant $C>0$ such that for any cube $Q$
\begin{align*}
\frac1{|Q|}\int_Qw(x)dx\le C\inf_{x\in Q}w(x).
\end{align*}
For $1<p<\infty$, we say that $w\in A_p$ if there exists a constant $C>0$ such that
\begin{align}\label{Ap}
[w]_{A_p}:=\sup_Q\bigg(\frac1{|Q|}\int_Qw(x)dx\bigg)\bigg(\frac1{|Q|}\int_Qw(x)^{1-p'}dx\bigg)^{p-1}\le C.
\end{align}
 We will adopt the following definition for the $A_\infty$ constant for a weight $w$ introduced by Fujii \cite{Fu} 
 and later by Wilson
 \cite{Wi}:
\begin{align}\label{Ainfty}
[w]_{A_\infty}:=\sup_Q\frac{1}{w(Q)}\int_QM(w\chi_Q)(x)\,dx.
\end{align}
Here $w(Q):=\int_Qw(x)\,dx,$ and the supremum above is taken over all cubes with edges parallel to the coordinate axes. When the supremum is finite, we
will say that $w$ belongs to the $A_\infty$ class.
 $A_\infty:=\bigcup_{p\ge1} A_p$. It is well known that if $w\in A_\infty$, then there exist $\delta\in(0,1]$ and $C>0$ such that for any interval $Q$ and
 measurable subset $E\subset Q$
\begin{align*}
\frac{w(E)}{w(Q)}\le C\bigg(\frac{|E|}{|Q|}\bigg)^\delta.
\end{align*}
Recall in Theorem \ref{thm5}, $$(w)_{A_p}:=\max\{[w]_{A_\infty},\,[w^{1-p'}]_{A_\infty}\}.$$ Using the facts that
$[w^{1-p'}]_{A_p'}^{1/p'}=[w]_{A_p}^{1/p}$ and $[w]_{A_\infty}\le C[w]_{A_p}$(see \cite{HP}), one easily checks that $$(w)_{A_p}\le
[w]_{A_p}^{\max(1,1/(p-1))}$$ (see \cite{HRT}).

\subsection{$\omega$-Dini Calder\'on--Zygmund   operators of convolution type}

A modulus of continuity is a function $\omega:[0,\infty)\rightarrow[0,\infty)$ with $\omega(0)=0$ that is subaddtive in the sense that
$$u\le t+s\Rightarrow \omega(u)\le \omega(t)+\omega(s).$$
Substituting $s=0$ one sees that $\omega(u)\le \omega(t)$ for all $0\le u\le t.$
Note that the composition and sum of two modulus of continuity are again  modulus of continuity.
In particular, if $\omega(t)$ is a modulus of continuity and $\theta\in (0,1),$ then $\omega(t)^\theta$ and $\omega(t^\theta)$ are also moduli of
continuity.
 The Dini norm of a modulus of continuity is defined by setting
 \begin{align} \|\omega\|_{Dini}:=\int_0^1\omega(t)\frac{dt}{t}<\infty.
\end{align} For any $c>0$ the integral can be equivalently (up to a $c$-dependent multiplicative constant) replaced by the sum over $2^{-j/c}$ with $j\in
\Bbb N.$
The basic example is $\omega(t)=t^\theta.$
 Let $T$ be a bounded linear operator on $L^2({\Bbb R}^n)$ represented as
  \begin{align} Tf(x)=\int_{{\Bbb R}^n}K(x-y)f(y)\,dy,\,\forall x \notin {\rm supp} f.
\end{align} We say that $T$ is an $\omega$-Dini Calder\'on--Zygmund  operators of convolution type if the kernel  $K$ satisfies the following size and
smoothness conditions:

for any $x\in {\Bbb R}^n\backslash \{0\}$,
\begin{align} \label{K0}|K(x)|\le \frac{C_K}{|x|^n};
\end{align}

for any $x,\,y\in {\Bbb R}^n$ with $2|y|\le |x|$,
\begin{align}\label{K1} |K(x-y)-K(x)|\le \frac{\omega(|y|/|x|)}{|x|^n}.
\end{align}

\quad

\subsection{Criterion for sharp weighted norm estimate}
We recall the criterion for sharp weighted norm estimate due to Lerner \cite{Ler 2}, which will be used frequently in the rest of the paper. The key role in this criterion is played by the grand maximal truncated operator $\mathcal{M}_U$, associated to a sub-linear operator $U$, defined by
\begin{align}\label{m}\mathcal{M}_{Uf}(x)=\mathcal{M}_Uf(x)=\sup_{Q\ni x}\, {\rm ess} \sup_{\xi\in Q} |U(f\chi_{{\Bbb R}^n\setminus{3Q}})(\xi)|,\end{align} where the supremum is
taken over all cubes $Q\subset {\Bbb R}^n$ containing $x.$
\begin{lemma}[\cite{Ler 0}]\label{lem 4} Assume that both $U$ and $\mathcal{M}_U$ be of weak type $(1,1)$, then for every compactly
supported $f\in L^p({\Bbb R}^n)\,(1<p<\infty)$ and $w\in A_p,$ \begin{align*} \|Uf\|_{L^p(w)}\lesssim K\{w\}_{A_p}\|f\|_{L^p(w)}, \end{align*} where $K=\|U\|_{L^1\rightarrow
L^{1,\infty}}+\|\mathcal{M}_U\|_{L^1\rightarrow L^{1,\infty}}.$
\end{lemma}

\subsection{The Stein--Weiss interpolation theorem with change of measures}
The following interpolation result with change of measures due to Stein and Weiss plays an important role in dealing with rough singular integrals. Again, we will use frequently this tool.
\begin{lemma}[\cite{sw}]\label{sw} Assume that $1\leq p_0,p_1\leq\infty$, that $w_0$ and $w_1$ be positive weights, and $T$ be a sublinear
operator  satisfying
$$
T:L^{p_i}(w_i)\rightarrow L^{p_i}(w_i),~~~~i=0,1,
$$
with quasi-norms $M_0$ and $M_1$, respectively. Then
$$
T:L^{p}(w)\rightarrow L^{p}(w),
$$
with quasi-norm $M\leq M_0^\lambda M_1^{1-\lambda}$, where
$$
\frac{1}{p}=\frac{\lambda}{p_0}+\frac{1-\lambda}{p_1},~~~w=w_0^{\frac{p\lambda}{p_0}}w_1^{\frac{p(1-\lambda)}{p_1}}.
$$\end{lemma}

\medskip
\section{Proof of Proposition \ref{pro2}}
We first give a useful Lemma which plays an important role in the proof of  Proposition \ref{pro2}.
\begin{lemma} \label{mw}   For any fixed $k\in \Bbb Z,$ let ${\Psi}_k$ be a function such that
{\rm supp} $ \Psi_k\subset\{x:|x| \lesssim 2^k\}.$ If ${\Psi}_k$ satisfy
 \begin{align}\label{psi0}\bigg\|(\sum_{k\in \Bbb Z}|{\Psi}_k\ast f|^2)^{1/2}\bigg\|_{L^{2}}\lesssim \gamma_1\|f\|_{L^2},\end{align}and
 \begin{align}\label{psi1}|{\Psi}_k(x)|\lesssim \frac{\gamma_2}{2^{kn}},\,\,\,k\in \Bbb Z,
\end{align}
 and for $0<2|y|\le |x|,$\begin{align}\label{psi2}\sum_{k\in \Bbb Z}|{\Psi}_k(x-y)-\Psi_k(x)|&\lesssim \gamma_3\frac{|y|^\gamma}{|x|^{n+\gamma}}.\end{align}
Then  for $1<p<\infty,$ and $w\in A_p,$
\begin{align}\label{m1} \bigg\|\bigg(\sum_{k\in \Bbb Z}|\Psi_k\ast f|^2\bigg)^{1/2}\bigg\|_{L^{p}(w)}&\lesssim  (\gamma_1+\gamma_2+\gamma_3)
\{w\}_{A_p}\|f\|_{L^{p}(w)}.\end{align}
\end{lemma}

\emph{Proof.}
To prove \eqref{m1}, by Lemma \ref{lem 4}, it suffices to verify that
\begin{align}\label{w1}\bigg\|\bigg(\sum_{k\in \Bbb Z}|{\Psi}_k\ast f|^2\bigg)^{1/2}\bigg\|_{L^{1,\infty}}&\lesssim
(\gamma_1+\gamma_2+\gamma_3)\|f\|_{L^1},\end{align} and
\begin{align}\label{w11} \|\mathcal{M}_{(\sum_{k\in \Bbb Z}|\Psi_k\ast f|^{2})^{1/2}}\|_{L^{1,\infty}}&\lesssim  (\gamma_1+\gamma_2+\gamma_3)
\|f\|_{L^{1}}.
\end{align}

To begin with, we verify \eqref{w1}.
Applying Calder\'{o}n--Zygmund  decomposition to $f$ at height  $\alpha$,  we obtain a disjoint  family
of dyadic cubes $\{Q\}$ with  total  measure
$$\dsum_Q|Q| \lesssim{\alpha}^{-1}\|f\|_{L^1},$$
which gives
$f=g+h,$
$\|g\|_{L^\infty}\lesssim\alpha,$
$\|g\|_{L^1}\lesssim\|f\|_{L^1},$ and
$h=\sum_{Q}h_Q,$
\text{supp}$(h_Q)\subseteq Q$,
$\int_{\mathbb R^n} h_Q(x)\,dx=0,$
$\sum\|h_Q\|_{L^1}\lesssim\|f\|_{L^1}.$
 By Chebychev's inequality and the $L^2$ estimate in \eqref{psi0} we get
 \begin{align*}
 \big| \big\{x: \big(\sum_{k\in \Bbb Z}|{\Psi}_k\ast g(x)|^2\big)^{1/2}\ge \alpha\big\}\big|\lesssim
 \frac{\big\|\big(\sum_{k\in \Bbb Z}|{\Psi}_k\ast g|^2\big)^{1/2}\big\|_{L^{2}}}{\alpha^2}\lesssim  \gamma_1\frac{\|g\|_{L^2}}{\alpha^2}\lesssim
 \gamma_1\frac{\|f\|_{L^1}}{\alpha}.
 \end{align*}
Then it suffices to estimate  $\big(\sum_{k\in \Bbb Z}|{\Psi}_k\ast h(x)|^2\big)^{1/2}$
away  from $\bigcup \widetilde{Q}$ with
 $\widetilde{Q}=2Q,$ where in general
we write $\rho Q$ for the $\rho$-dilated $Q$, that is, for the cube with same center as $Q$ and with side
length $\widetilde{Q}$.  Thus  by \eqref{psi2}
\begin{align*}&\alpha\big|\{x\notin\bigcup \widetilde{Q}:\big(\sum_{k\in \Bbb Z}|{\Psi}_k\ast h(x)|^2\big)^{1/2}>\alpha\}\big|\\
&\le\dsum_{Q}\dint_{x\notin \widetilde{Q}}\dint_Q|h_Q(y)|\sum_{k\in \Bbb Z}\big|{\Psi}_k(x-y)-{\Psi}_k(x-y_Q)\big|dy\,dx\\
&\lesssim\dsum_{Q}\dint_Q|h_Q(y)|\dint_{|x-y_Q|\ge 2\ell(Q)}\frac{\gamma_3|y-y_Q|^\gamma}{|x-y|^{n+\gamma}}\,dxdy\\
&\lesssim \gamma_3\dsum_{Q}\|h_Q\|_{L^1}\\
&\lesssim   \gamma_3\|f\|_{L^1},\end{align*}where $y_Q$ denotes the center of $Q.$ Thus we finish the proof of
\eqref{w1}.

\quad

Next, we verify \eqref{w11}. Take $x$ and $\xi $ in any fixed cube $Q.$ Let $B(x)=B(x,2\sqrt{n}\ell(Q)),$ then $3Q\subset B_x.$
By using the triangle  inequality, we get \begin{align*}
& \big(\sum_{k\in \Bbb Z}|\Psi_k\ast (f{\chi_{{\Bbb R}^n\setminus 3Q}})(\xi)|^2\big)^{1/2}\\& \le \big|\big(\sum_{k\in \Bbb Z}|\Psi_k\ast (f{\chi_{{\Bbb
R}^n\setminus B_x}})(\xi)|^2\big)^{1/2}-\big(\sum_{k\in \Bbb Z}|\Psi_k\ast (f{\chi_{{\Bbb R}^n\setminus B_x}})(x)|^2\big)^{1/2}\big|+\big(\sum_{k\in \Bbb
Z}|\Psi_k\ast f{\chi_{B_x\setminus 3Q}})(\xi)|^2\big)^{1/2}\\&\quad+\big(\sum_{k\in \Bbb Z}|\Psi_k\ast (f{\chi_{{\Bbb R}^n\setminus
B_x}})(x)|^2\big)^{1/2}\\&=:I+I\!I+I\!I\!I.
\end{align*}
 For the term $I,$  by \eqref{psi2} we have
\begin{align*}
I&\le\int_{{\Bbb R}^n \setminus B_x} \sum_{k\in \Bbb Z}|\Psi_k(\xi-y)-\Psi_k(x-y)||f(y)|\,dy\le\int_{{\Bbb R}^n \setminus
B_x}\frac{\gamma_3|x-\xi|}{|x-y|^{n+1}}|f(y)|\,dy
\lesssim \gamma_3Mf(x).\nonumber
\end{align*}
For the term $I\!I$,  by the size condition \eqref{psi1} and the fact that $|x-y|\simeq |\xi-y|$ (since $2|x-\xi|\le |x-y|$), we get
\begin{align*}
I\!I&\le \int_{{\Bbb R}^n}\sum_{k\in \Bbb Z}|\Psi_k(\xi-y)||f|{\chi_{B_x\setminus 3Q}}(y)\,dy\lesssim  \gamma_2\int_{3\ell(Q)\le|x-y|\le
2\sqrt{n}\ell(Q)}|\xi-y|^{-n}|f(y)|\,dy\lesssim \gamma_2Mf(x).
\end{align*}
For the last term, recalling that supp $\Psi_k\subset \{x:|x|\lesssim 2^{k}\}$ and that $B_{x}=\{y:|y-x|\leq2\sqrt{n} \ell(Q)\}$,
we have that
\begin{align*}
\Psi_k*f\chi_{\mathbb{R}^{n}\setminus B_{x}}&=\int_{2\sqrt{n} \ell(Q)<|x-y|\lesssim 2^{k}}\Psi_k(x-y)f(y)dy.
\end{align*}
Thus,
\begin{align*}
I\!I\!I
&= \bigg(\sum_{2\sqrt{n} \ell(Q)\lesssim2^{k}}|\Psi_k\ast f\chi_{\mathbb{R}^{n}\setminus B_{x}}(x)|^{2}\bigg)^{1/2}\leq \bigg(\sum_{k\in \Bbb Z}|\Psi_k\ast
f(x)|^{2}\bigg)^{1/2}+\bigg(\sum_{2\sqrt{n}  \ell(Q)\lesssim2^{k}}|\Psi_k\ast f\chi_{ B_{x}}(x)|^{2}\bigg)^{1/2}.
\end{align*}
By using the size estimate \eqref{psi1}, we have
\begin{align*}
\bigg(\sum_{2\sqrt{n}  \ell(Q)\lesssim2^{k}}|\Psi_k\ast f\chi_{ B_{x}}(x)|^{2}\bigg)^{1/2}
&\lesssim\sum_{2\sqrt{n} \ell(Q)\lesssim2^{k}}\frac{\gamma_3}{2^{kn}}\int_{|x-y|\le2\sqrt{n} \ell(Q)}|f(y)|dy\\
&\lesssim \frac{\gamma_2}{\ell(Q)^{n}}\int_{|x-y|\le2\sqrt{n} \ell(Q)}|f(y)|dy\\
&\lesssim \gamma_2Mf(x),
\end{align*}
which gives
$$I\!I\!I
\leq \bigg(\sum_{k\in \Bbb Z}|\Psi_k\ast f(x)|^{2}\bigg)^{1/2}+\gamma_2Mf(x).$$
Combining the estimates of $I,\,I\!I$ and $I\!I\!I,$ we get
\begin{align*} \mathcal{M}_{(\sum_{k\in \Bbb Z}|\Psi_k\ast f(x)|^{2})^{1/2}}&\lesssim  \bigg(\sum_{k\in \Bbb Z}|\Psi_k\ast
f(x)|^{2}\bigg)^{1/2}+(\gamma_2+\gamma_3)Mf(x).
\end{align*}
Then the weak type $(1,1)$ of the Hardy--Littlewood maximal function $M$ and the estimate in  \eqref{w1}  imply
\begin{align*} \|\mathcal{M}_{(\sum_{k\in \Bbb Z}|\Psi_k\ast f|^{2})^{1/2}}\|_{L^{1,\infty}}&\lesssim  (\gamma_1+\gamma_2+\gamma_3) \|f\|_{L^{1}},
\end{align*}
which shows that \eqref{w11} holds.

The proof of Lemma \ref{mw} is complete. \qed

\quad

Now we return to the proof of Proposition \ref{pro2}.
We point out that to estimate $\{\phi_{k}*K_{k}\},$ we do not need any regularity condition on $K,$
we only assume $K$ satisfies the size condition \eqref{K0} and the cancellation condition \eqref{2}.

We  first verify the size estimate of $\phi_k\ast
K_k(x)$ for any fixed $k\in \Bbb Z.$ Since supp  ${\phi}_k\subset \{x: |x|\le 2^k/4\},$ we have supp $\phi_k\ast K_k\subset\{x:|x|\le 2^{k+1}\}.$
  Then from \eqref{K0} and \eqref{2} we have
\begin{align}\label{2j}|{{\phi}}_k\ast K_k(x)|&=\bigg|\int_{|z|\le 2^{k}}{\phi}_k(x-z)K(z)\,dz\bigg|\\
&\le
\bigg|\int_{|z|\le2^{k}} ({\phi}_k(x-z)-{\phi}_k(x))K(z)\,dz\bigg|\nonumber\\
&\lesssim \frac{C_K}{2^k}\int_{|z|\le2^{k}} |(\nabla\phi)_k(x-\theta
z)|\frac{1}{|z|^{n-1}}\,dz
\nonumber\\
&\lesssim\frac{C_K}{2^{k(n+1)}}\int_{|z|\le2^{k}} \frac{1}{|z|^{n-1}}\,dz\nonumber\\&\lesssim\frac{C_K}{2^{k{n}}}\chi_{|x|\le 2^{k+1}}(x),\nonumber
\end{align}
which shows that  \eqref{sd} holds.

Now, we estimate $\nabla\phi_{k}\ast K_k(x)$.
Recall that $\phi\in C_0^\infty({\Bbb R}^n)$  with $\int \phi =1$ and supp $\phi\subset \{x: |x|\le 1/4\}.$ It is easy to verify that $\int\nabla \phi=0$
and supp $(\nabla \phi)\subset \{x: |x|\le 1/4\}$ with $\nabla \phi\in C_0^\infty({\Bbb R}^n).$ Then by writing
\begin{align*}\nabla\phi_{k}\ast K_k(x)&=\frac{1}{2^k}(\nabla\phi)_{k}\ast K_k(x)
\end{align*}
and repeating the argument of \eqref{2j}, we get
\begin{align}\label{lj1}|\nabla\phi_{k}\ast K_k(x)|&=\frac{1}{2^k}|(\nabla\phi)_{k}\ast K_k(x)|\lesssim\frac{C_K}{2^{k{(n+1)}}}\chi_{|x|\le 2^{k+1}}(x).
\end{align}
Thus
for $|x|\ge 2|y|,$
 \begin{align*} \sum_{k\in \Bbb Z}|{{\phi}}_{k}\ast K_k(x-y)-{{\phi}}_{k}\ast K_k(x)|\lesssim\sum_{k\in \Bbb Z}|\nabla{\phi}_k\ast K_k(x-\theta y
 )|y|\lesssim C_K\frac{|y|}{|x|^{n+1}},
\end{align*} which shows that  \eqref{smd} holds.

On the other hand, it is easy to verify that \begin{align*}|{{\phi}}_k\ast K_k(x)|&\lesssim\frac{C_K}{2^{k{n}}}\chi_{|x|\le 2^{k+1}}\lesssim C_K\frac{2^k}{(2^k+|x|)^{n+1}},
\end{align*}
and that
\begin{align*}|{{\phi}}_k\ast K_k(x+h)-{{\phi}}_k\ast K_k(x)|&\lesssim C_K\frac{|h|^\theta}{(2^k+|x|)^{n+\theta}}, \qquad |h|\le 2^k
\end{align*} for some $\theta\in (0,1)$.
In addition, since ${{\phi}}_k\ast K_k\ast 1=0,$ by applying the $T1$ Theorem for square functions (see \cite{CHJ}) we have
\begin{align}\label{w2}\bigg\|\bigg(\sum_{k\in \Bbb Z}|\phi_k\ast K_k\ast f|^2\bigg)^{1/2}\bigg\|_{L^{2}}&\lesssim C_K\|f\|_{L^2},\end{align} which shows that  \eqref{s2d} holds.

Then  applying  Lemma \ref{mw} to $\{\Psi_k\}=\{\phi_k\ast K_k\},$ we get
\begin{align*} \bigg\|\bigg(\sum_{k\in \Bbb Z}|\phi_{k}\ast K_k\ast f|^2\bigg)^{1/2}\bigg\|_{L^{1,\infty}}&\lesssim
C_K\|f\|_{L^{1}},\end{align*}\begin{align*} \bigg\|\mathcal{M}_{(\sum_{k\in \Bbb Z}|\phi_{k}\ast K_k\ast f|^2)^{1/2}}\bigg\|_{L^{1,\infty}}&\lesssim
C_K\|f\|_{L^{1}}\end{align*}and for $1<p<\infty$ and $w\in A_p,$\begin{align*} \bigg\|\bigg(\sum_{k\in \Bbb Z}|\phi_{k}\ast K_k\ast f|^2\bigg)^{1/2}\bigg\|_{L^{p}(w)}&\lesssim
C_K\{w\}_{A_p}\|f\|_{L^{p}(w)}.\end{align*}This verifies \eqref{w1d}, \eqref{w1md} and \eqref{swd}.

The proof of  Proposition \ref{pro2} is complete.\qed

\quad

\section{Proof of Proposition \ref{pro4}}

\textbf{Step 1.} First, we would like to prove that  \begin{align}\label{w3}\bigg\|\bigg(\sum_{k\in \Bbb Z}|({{\phi}}_k-\delta_0)\ast K^k\ast
f|^q\bigg)^{1/q}\bigg\|_{L^{2}}&\lesssim(C_{K}+\|\omega\|_{Dini}+\|\omega\|_{Dini}^{1/q'}\|\omega^{1/2}\|_{Dini}^{2/q})\|f\|_{L^2}.\end{align}   We  first use Fourier transform and Plancherel theorem to deal
with the $L^2$-norm of $(\sum_{k\in \Bbb Z}|({{\phi}}_k-\delta_0)\ast K^k\ast f|^2)^{1/2}.$  Since $K$ is not necessary of homogeneous type, the computation
of $|\widehat{K^k}(\xi)|$ is not trivial. Write
\begin{align*}&\widehat{K^k}(\xi)=\int_{2^{k}<|x|<|\xi|^{-1}}K(x)e^{-2\pi ix \cdot \xi} dx
+\int_{|\xi|^{-1}<|x|}
K(x) e^{-2\pi i x \cdot \xi} dx
=: I_{1}(\xi)+I_{2}(\xi).
\end{align*}
For the term $I_1$, by the size estimate \eqref{K0} and cancellation \eqref{2}, we have
\begin{align*}
|I_{1}(\xi)|&=
\bigg|\int_{2^{k}<|x|\le|\xi|^{-1}}K(x)(e^{-2\pi ix \cdot \xi}-1) dx\bigg|\\
&\lesssim2\pi|\xi|\int_{|x|\le|\xi|^{-1}}|x||K(x)|dx\\
&\lesssim C_{K}
\end{align*}
uniformly in $k$.

We now consider $I_{2}(\xi)$.
Let $z=\frac{\xi}{2|\xi|^2}$ so that $e^{2\pi i z \cdot \xi}=-1$
and $2|z|=|\xi|^{-1}$.
By changing variables $x=x'-z$, we rewrite $I_2(\xi)$ as
\begin{align*}
I_{2}(\xi)=-\int_{|\xi|^{-1}<|x'-z|}K(x'-z)e^{-2\pi ix'\cdot\xi}dx'.
\end{align*}
Taking the average of the above inequality and the original definition of $I_2$ gives
\begin{align*}
I_{2}(\xi)=\frac{1}{2}\int_{|\xi|^{-1}<|x|}K(x)e^{-2\pi i x \cdot \xi}dx
-\frac{1}{2}\int_{|\xi|^{-1}<|x-z|}K(x-z)e^{-2\pi i x \cdot \xi} dx.
\end{align*}
Now we split
\begin{align*}
I_{2}(\xi)&=
\frac{1}{2}\int_{|\xi|^{-1}<|x|}(K(x)-K(x-z))e^{-2\pi i x\cdot\xi} dx\\
&+\frac{1}{2}\int_{{\Bbb R}^n} K(x-z)(\chi_{\{|\xi|^{-1}<|x|\}}-\chi_{\{|\xi|^{-1}<|x-z|\}})
e^{-2\pi i x\cdot \xi}dx\\
&=: J_{1}(\xi)+J_{2}(\xi).
\end{align*}
For the term $J_1$, by noting that  $2|z|=|\xi|^{-1}$ and using \eqref{K1} we get
\begin{align*}
|J_{1}(\xi)|&\leq\frac{1}{2}\int_{|\xi|^{-1}<|x|}\omega\left(\frac{|z|}{|x|}\right)dx
\leq \frac{1}{2}\int_{|\xi|^{-1}<|x|}\omega\left(\frac{|\xi|^{-1}}{2|x|}\right)dx
\lesssim\|\omega\|_{Dini}.
\end{align*}
For the term $J_{2}$, we note that
$\chi_{|\xi|^{-1}<|x|}-\chi_{|\xi|^{-1}<|x-z|}$ is nonzero if and
only if $|\xi|^{-1}-|z|\le|x-z|\le|\xi|^{-1}+|z|$.
Thus by \eqref{K0}
\begin{align*}
|J_{2}(\xi)|&\leq
  \frac{1}{2}\int_{|\xi|^{-1}-|z|\le|x-z|\le|\xi|^{-1}+|z|} |K(x-z)|dx\\
&\leq
\frac{1}{2}\int_{|\xi|^{-1}-|z|\le|x|\le|\xi|^{-1}+|z|}\frac{C_{K}}{|x|^{n}}dx \\
&\leq
\frac{1}{2}\int_{\frac{|\xi|^{-1}}{2}\le|x|\le\frac{3|\xi|^{-1}}{2}}\frac{C_{K}}{|x|^{n}}dx\\
&\lesssim  C_{K}.
\end{align*}
Combining the estimates of $J_{1}(\xi)$ and $J_{2}(\xi),$
we get
\begin{align*}
|I_{2}(\xi)|\lesssim C_{K}+\|\omega\|_{Dini},
\end{align*}
which, together with the estimate of  $I_{1}(\xi)$, gives
\begin{align}\label{kj}
|\widehat{K^{k}}(\xi)|\lesssim C_{K}+\|\omega\|_{Dini}.
\end{align}

 For $|2^k\xi|>1,$ write
\begin{align*}&\widehat{K^k}(\xi)=\int_{2^k<|x|}
K(x) e^{-2\pi i x \cdot \xi} dx.
\end{align*}
Let $z=\frac{\xi}{2|\xi|^2}$ so that $e^{2\pi i z \cdot \xi}=-1$
and $2|z|=|\xi|^{-1}$. Similarly to $I_2(\xi)$, we get
\begin{align*}
\widehat{K^k}(\xi)&=
\frac{1}{2}\int_{2^k<|x|}(K(x)-K(x-z))e^{-2\pi i x\cdot\xi} dx\\
&+\frac{1}{2}\int_{{\Bbb R}^n} K(x-z)(\chi_{\{2^k<|x|\}}-\chi_{\{2^k<|x-z|\}})
e^{-2\pi i x\cdot \xi}dx\\
&=: J_{1,k}(\xi)+J_{2,k}(\xi).
\end{align*}
For the term $J_{1,k}$, by noting that  $2|z|=|\xi|^{-1}$ and using \eqref{K1} we get
\begin{align*}
|J_{1,k}(\xi)|&\leq\frac{1}{2}\int_{2^k<|x|}\omega\left(\frac{|z|}{|x|}\right)dx
\leq \frac{1}{2}\int_{2^k<|x|}\omega\left(\frac{|\xi|^{-1}}{2|x|}\right)dx
\lesssim\|\omega^{1/2}\|_{Dini}\omega^{1/2}(|2^k\xi|^{-1}).
\end{align*}
For the term $J_{2,k}$, we note that
$\chi_{2^k<|x|}-\chi_{2^k<|x-z|}$ is nonzero if and
only if $2^k-|z|\le|x-z|\le2^k+|z|$.
Thus by \eqref{K0}
\begin{align*}
|J_{2,k}(\xi)|&\leq
  \frac{1}{2}\int_{2^k-|z|\le|x-z|\le2^k+|z|} |K(x-z)|dx\\
&\leq
\frac{1}{2}\int_{2^k-|z|\le|x|\le2^k+|z|}\frac{C_{K}}{|x|^{n}}dx \\
&\leq
\frac{1}{2^k}\int_{2^k-|z|\le|x|\le2^k+|z|}\frac{C_{K}}{|x|^{n-1}}dx\\
&\lesssim  C_{K}|2^k\xi|^{-1}.
\end{align*}
Combining the estimates of $J_{1,k}(\xi)$ and $J_{2,k}(\xi),$
we get
\begin{align}\label{Kk}
|\widehat{K^k}(\xi)|\lesssim C_{K}|2^k\xi|^{-1}+\|\omega^{1/2}\|_{Dini}\omega^{1/2}(|2^k\xi|^{-1}).
\end{align}
By the Plancherel Theorem,  we get  \begin{align*}
\bigg\|\bigg(\sum_{k\in \Bbb Z}|(\delta_0-\phi_k)\ast K^k\ast f|^2\bigg)^{1/2}\bigg\|_{ L^2}^2
&=\sum_{k\in \mathbb{Z}}\int_{{\Bbb R}^n}|\widehat{(\delta_0-\phi_k)}(\xi)|^2|\widehat{K^k}(\xi)|^2|\widehat{f}(\xi)|^2d\xi.
\end{align*}
Since $\widehat{\phi}\in \mathcal{S}({\Bbb R}^n)$ and $\widehat{\phi}(0)=1,$ then $|1-\widehat{\phi_k}(\xi)|\lesssim \min(|2^k\xi|, 1).$  Together with   \eqref{kj} and \eqref{Kk}, we have
\begin{align*}
&\sum_{k\in \mathbb{Z}}|\widehat{(\delta_0-\phi_k)}(\xi)|^2|\widehat{K^k}(\xi)|^2|\widehat{f}(\xi)|^2\\
&\lesssim \big((C_{K}+\|\omega\|_{Dini})^2\sum_{|2^k\xi|\leq1}|2^k\xi|^{2}+C_K^2\sum_{|2^k\xi|\geq1}|2^{k}\xi|^{-2}+\|\omega^{1/2}\|_{Dini}^2\sum_{|2^k\xi|\geq1}\omega(|2^k\xi|^{-1})\big)|\widehat{f}(\xi)|^2\\
&\lesssim((C_{K}+\|\omega\|_{Dini})^2+C_K^2+\|\omega^{1/2}\|_{Dini}^2\|\omega\|_{Dini})|\widehat{f}(\xi)|^2,
\end{align*}
which implies that
\begin{align}\label{kk2}
\bigg\|\bigg(\sum_{k\in \Bbb Z}|(\delta_0-\phi_k)\ast K^k\ast f|^2\bigg)^{1/2}\bigg\|_{ L^2}\lesssim (C_{K}+\|\omega\|_{Dini}+\|\omega^{1/2}\|_{Dini}\|\omega\|_{Dini}^{1/2})\|f\|_{L^2}.
\end{align}

On the other hand, we claim that   for any fixed $k\in \Bbb Z,$
\begin{align}\label{uj00}|(\delta_0-\phi_k)\ast K^k(x)|&\lesssim\omega(2^k/|x|)\frac{1}{|x|^{n}}\chi_{|x|>3\cdot2^k/4}+\frac{C_K}{|x|^n}\chi_{\frac{3}{4}\cdot2^k\le |x|\le \frac{5}{4}\cdot2^k},
\end{align} which will be proved in Step 2. Then we can  get
\begin{align}\label{kk3}
\big\|\sup_{k\in \Bbb Z}|(\delta_0-\phi_k)\ast K^k\ast f|\big\|_{ L^2}\lesssim(C_{K}+\|\omega\|_{Dini})\|Mf\|_{L^2}\lesssim (C_{K}+\|\omega\|_{Dini})\|f\|_{L^2}.
\end{align}
Interpolating between \eqref{kk2} and \eqref{kk3}, we get for $q>2$
\begin{align*}
\bigg\|\bigg(\sum_{k\in \Bbb Z}|(\delta_0-\phi_k)\ast K^k\ast f|^q\bigg)^{1/q}\bigg\|_{ L^2}\lesssim (C_{K}+\|\omega\|_{Dini}+\|\omega\|_{Dini}^{1/q'}\|\omega^{1/2}\|_{Dini}^{2/q})\|f\|_{L^2}.
\end{align*} we get \eqref{w3}, which shows \eqref{s2u}.

\textbf{Step 2.} We  now estimate $|(\delta_0-\phi_k)\ast K^k(x)|$ for any fixed $k\in \Bbb Z.$
\begin{align*}|(\delta_0-\phi_k)\ast K^k(x)|&=\bigg|\int_{{\Bbb R}^n} {\phi}_k(x-z)(K^k(z)-K^k(x))\,dz\bigg|\\&\le\bigg|\int_{{\Bbb R}^n} {\phi}_k(x-z)(K(z)-K(x))\chi_{|z|>
2^{k}}\,dz\bigg|+\bigg|\int_{{\Bbb R}^n} {\phi}_k(x-z)K(x))(\chi_{|z|> 2^{k}}-\chi_{|x|> 2^{k}})\,dz\bigg|\\&=:I_{k,1}(x)+I_{k,2}(x).
\end{align*}
We first consider $I_{k,1}$. By noting that supp $ {\phi}_k\subset \{x: |x|\le 2^k/4\}$  and that $|z|>2^{k},$ we have $|x|\ge \frac{3|z|}{4}> 3\cdot2^k/4.$ Then by the regularity condition of $K$ in
\eqref{K1},\begin{align}\label{2jj}I_{k,1}(x)&\lesssim\int_{|z|>2^{k}} |{\phi}_k(x-z)|\frac{\omega(|
x-z|/|x|)}{|x|^{n}}\,dz\lesssim\frac{\omega(2^k/|x|)}{|x|^{n}}\chi_{|x|>3\cdot2^k/4}(x).
\end{align}
For the term $I_{k,2}$, we first note that  $\chi_{|z|> 2^{k}}-\chi_{|x|> 2^{k}}$ is nonzero if and only if $\frac{3}{4}\cdot2^k\le |x|\le \frac{5}{4}\cdot2^k$ since $|x-z|\le 2^k/4$. Thus,
\begin{align*}I_{k,2}(x)\lesssim
|K(x)|\chi_{\frac{3}{4}\cdot2^k\le |x|\le \frac{5}{4}\cdot2^k}(x)\int_{{\Bbb R}^n} |{\phi}_k(x-z)|\,dz\lesssim \frac{C_K}{|x|^n}\chi_{\frac{3}{4}\cdot2^k\le |x|\le \frac{5}{4}\cdot2^k}(x).
\end{align*}
Combining the two cases,  we get that for any fixed $k\in \Bbb Z,$
\begin{align}\label{uj0}|({{\phi}}_k-\delta_0)\ast K^k(x)|&\lesssim\omega(2^k/|x|)\frac{1}{|x|^{n}}\chi_{|x|>3\cdot2^k/4}+\frac{C_K}{|x|^n}\chi_{\frac{3}{4}\cdot2^k\le |x|\le \frac{5}{4}\cdot2^k},
\end{align}
which verifies \eqref{uj00} and implies that
\begin{align}\label{tk1}\sum_{k\in \Bbb Z}|({{\phi}}_k-\delta_0)\ast K^k(x)|&\lesssim (\|\omega\|_{Dini}+C_K)\frac{1}{|x|^{n}}.
\end{align}
Hence \eqref{su} holds.

We now estimate $\|\{({{\phi}}_k-\delta_0)\ast K^k(x-y)-({{\phi}}_k-\delta_0)\ast K^k(x)\}_{k}\|_{\ell^q}$ for $0<|y|\le |x|/2.$ To begin with,
we claim that for $|y|\le \frac{2^k}{2},$
\begin{align}\label{j1}
&|({{\phi}}_k-\delta_0)\ast K^k(x-y)-({{\phi}}_k-\delta_0)\ast K^k(x)|\\&\lesssim\bigg(\frac{\omega(|y|/|x|)}{|x|^{n}}+
C_K\frac{|y|^\theta}{|x|^{n+\theta}}\bigg)\chi_{|x|\ge2^{k}/2}(x)+\frac{C_K}{|x|^n}\chi_{2^{k}-|y|\leq|x|\le2^{k}+|y|}(x).\nonumber
\end{align}
Assuming \eqref{j1} for the moment,
we write
\begin{align*}&\bigg(\sum_{k\in \Bbb Z}|({{\phi}}_{k}-\delta_0)\ast K^k(x-y)-({{\phi}}_{k}-\delta_0)\ast K^k(x)|^q\bigg)^{1/q}\\&\le\bigg(\big(\sum_{
{\frac{2^k}{2}}\le |y|}+\sum_{ {\frac{2^k}{2}}\ge |y|}\big)|({{\phi}}_{k}-\delta_0)\ast K^k(x-y)-({{\phi}}_{k}-\delta_0)\ast K^k(x)|^q\bigg)^{1/q}\\&=:I\!I_1+I\!I_2.
\end{align*}
For the term $I\!I_1$, by noting that $|x|\geq 2|y|>0$ and by using \eqref{uj0} we get
\begin{align*}I\!I_1&\le\bigg(\sum_{ \frac{{2^k}}{2}\le |y|}\big(|({{\phi}}_{k}-\delta_0)\ast K^k(x)|+|({{\phi}}_{k}-\delta_0)\ast
K^k(x-y)|\big)^q\bigg)^{1/q}\\&\lesssim\sum_{ \frac{{2^k}}{2}\le|y|}\frac{C_K}{|x|^n}\chi_{\frac{3}{4}\cdot2^k\le |x|\le \frac{5}{4}\cdot2^k}(x)+\bigg(\sum_{ \frac{{2^k}}{2}\le
|y|}\big(\frac{\omega(2^k/|x|)}{|x|^{n}}\chi_{|x|>3\cdot2^k/4}\big)^{q}\bigg)^{1/q}\\&\lesssim
\frac{C_K|y|}{|x|^{n+1}}+\frac{\omega(|y|/|x|)^{\frac{1}{q'}}}{|x|^{n}}\bigg(\sum_{|x|>3\cdot2^k/4}\omega(2^k/|x|)\bigg)^{1/q}\\&\lesssim
\frac{C_K|y|}{|x|^{n+1}}+\frac{\omega(|y|/|x|)^{\frac{1}{q'}}}{|x|^{n}}\|\omega\|_{Dini}^{1/q}.
\end{align*}
For the term $I\!I_2$, by using \eqref{j1},
\begin{align*}I\!I_2&\lesssim\sum_{\frac{{2^k}}{2}\ge  |y|}C_K\frac{|y|^\theta}{|x|^{n+\theta}}\chi_{|x|\ge 2^{k}/4}+\bigg(\sum_{\frac{{2^k}}{2}\ge
|y|}\bigg(\frac{\omega(|y|/|x|)}{|x|^{n}}\chi_{|x|\ge\frac{2^{k}}{4}}\bigg)^q\bigg)^{1/q}+\sum_{\frac{{2^k}}{2}\ge
|y|}\frac{C_K}{|x|^n}\chi_{2^{k}-|y|\le|x|\le2^{k}+|y|}\\&\lesssim
\frac{C_K |y|^{\theta/2}}{|x|^{n+\theta/2}}+\frac{\omega(|y|/|x|)^{\frac{1}{q'}}}{|x|^{n}}\bigg(\sum_{{\frac{2^k}{2}}\ge
|y|}\omega(|y|/2^k)\bigg)^{1/q}+\sum_{{\frac{2^k}{2}}\ge  |y|}\frac{C_K}{|x|^n}\chi_{2^{k}-|y|\le|x|\le2^{k}+|y|}\\&\lesssim\frac{C_K
|y|^{\theta/2}}{|x|^{n+\theta/2}}+\frac{\omega(|y|/|x|)^{\frac{1}{q'}}}{|x|^{n}}\|\omega\|_{Dini}^{1/q}+\sum_{{\frac{2^k}{2}}\ge
|y|}\frac{C_K}{|x|^n}\chi_{2^{k}-|y|\le|x|\le2^{k}+|y|}.
\end{align*}Combining  the estimates of $I\!I_1$ and $I\!I_2$,  we get
 \begin{align}\label{jkx}&\bigg(\sum_{k\in \Bbb Z}|({{\phi}}_{k}-\delta_0)\ast K^k(x-y)-({{\phi}}_{k}-\delta_0)\ast K^k(x)|^q\bigg)^{1/q}\\&\lesssim\frac{C_K
 |y|^{\theta/2}}{|x|^{n+\theta/2}}+\frac{\omega(|y|/|x|)^{\frac{1}{q'}}}{|x|^{n}}\|\omega\|_{Dini}^{1/q}+\sum_{\frac{{2^k}}{2}\ge
 |y|}\frac{C_K}{|x|^n}\chi_{2^{k}-|y|\le|x|\le2^{k}+|y|},\nonumber
\end{align}which verifies \eqref{smu}.

\quad

Now we return to give the proof of \eqref{j1}. We write
\begin{align*}
&|({{\phi}}_k-\delta_0)\ast K^k(x-y)-({{\phi}}_k-\delta_0)\ast K^k(x)|\\
&=\bigg|\int_{{\Bbb R}^n}\phi_{k}(x-z)\big(K^{k}(z-y)-K^{k}(x-y)\big)\,dz-\int_{{\Bbb R}^n}\phi_{k}(x-z)\big(K^{k}(z)-K^{k}(x)\big)\,dz\bigg|\\
&\leq\bigg|\int_{{\Bbb R}^n}\phi_{k}(x-z)\big(K^{k}(x-y)-K^{k}(x)\big)\,dz\bigg|+\bigg|\int_{{\Bbb R}^n}\phi_{k}(x-z)\big(K^{k}(z-y)-K^{k}(z)\big)\,dz\bigg|\\
&=:I\!I\!I_1+I\!I\!I_2.
\end{align*}
For the term $I\!I\!I_1$, recall that  $|y|\le\frac{2^k}{2}$. By noting the facts that $0<|y|\le\frac{2^k}{2}$ and $|x-y|\ge 2^{k}$ imply $|x|>\frac{2^k}{2}$ and that $\chi_{|x-y|>2^{k}}-\chi_{|x|>2^{k}}$ is nonzero  if and only if
 $2^{k}-|y|\leq|x|\le 2^{k}+|y|$, we have
\begin{align*}
I\!I\!I_1
&\lesssim\int_{{\Bbb R}^n}|\phi_{k}(x-z)||K(x-y)-K(x)|\chi_{|x-y|>2^{k}}\,dz\\
&\quad+\int_{{\Bbb R}^n}|\phi_{k}(x-z)||K(x)||\chi_{|x-y|>2^{k}}-\chi_{|x|>2^{k}}|\,dz\\
&\lesssim \frac{\omega(|y|/|x|)}{|x|^{n}}\chi_{|x|>\frac{2^{k}}{2}}+\frac{C_K}{|x|^n}|\chi_{|x-y|>2^{k}}-\chi_{|x|>2^{k}}|\\
&\lesssim \frac{\omega(|y|/|x|)}{|x|^{n}}\chi_{|x|>\frac{2^{k}}{2}}+\frac{C_K}{|x|^n}\chi_{2^{k}-|y|\leq|x|\le2^{k}+|y|}.\end{align*}
Similarly, for the term $I\!I\!I_2$, we get
\begin{align*}
I\!I\!I_2
&\lesssim\int_{{\Bbb R}^n}|\phi_{k}(x-z)|\frac{\omega(|y|/|z|)}{|z|^{n}}\chi_{|z|>\frac{2^{k}}{2}}\,dz+C_K\int_{{\Bbb R}^n}\frac{|\phi_{k}(x-z)|}{|z|^n}\chi_{2^{k}-|y|\le|z|\leq2^{k}+|y|}\,dz.
\end{align*}
Since $|x-z|\le2^k/4$ and $|z|>\frac{2^k}{2}$, we get
$
3|z|/2\ge|x|\ge |z|/2>2^k/4
$. Thus, there exists some $\theta\in (0,1)$ and $\frac{1}{p'}=\theta,$
\begin{align*}
I\!I\!I_2&\lesssim \frac{\omega(|y|/|x|)}{|x|^{n}}\chi_{|x|\ge2^{k}/4}+\frac{C_K}{|x|^{n+\theta}}\chi_{|x|\ge2^{k}/4} \bigg(\int_{{\Bbb R}^n}|\phi_k(x-z)|^p\,dz\bigg)^{\frac{1}{p}}\bigg(\int_{2^{k}-|y|\leq|z|\le2^{k}+|y|}|z|^{\theta p'}\,dz\bigg)^{\frac{1}{p'}}\\
&\lesssim\frac{\omega(|y|/|x|)}{|x|^{n}}\chi_{|x|\ge2^{k}/4}+C_K\frac{|y|^\theta}{|x|^{n+\theta}}\chi_{|x|\ge2^{k}/4}.
\end{align*}
Combining the estimates of $I\!I\!I_1$ and $I\!I\!I_2$, we get \eqref{j1}.

\quad\quad

\textbf{Step 3.} By Lemma \ref{lem 4},
to prove \eqref{swu}, we need to verify that  \begin{align}\label{wj1} \big\|(\sum_{k\in \Bbb Z}|(\phi_{k}-\delta_0)\ast K^k\ast
f|^q)^{1/q}\big\|_{L^{1,\infty}}\lesssim  (C_K+\|\omega\|_{Dini}+\|\omega\|_{Dini}^{1/q}\|\omega^{1/q'}\|_{Dini}+\|\omega\|_{Dini}^{1/q'}\|\omega^{1/2}\|_{Dini}^{2/q})\|f\|_{L^{1}}.\end{align}and\begin{align}\label{wj2}
\|\mathcal{M}_{(\sum_{k\in \Bbb Z}|(\phi_{k}-\delta_0)\ast K^k\ast f|^q)^{1/q}}\|_{L^{1,\infty}}\lesssim
(C_K+\|\omega\|_{Dini}+\|\omega\|_{Dini}^{1/q}\|\omega^{1/q'}\|_{Dini}+\|\omega\|_{Dini}^{1/q'}\|\omega^{1/2}\|_{Dini}^{2/q})\|f\|_{L^{1}}.\end{align}
The two above inequalities just are \eqref{w1u} and \eqref{w1mu}.

To verify \eqref{wj1}, we apply
the Calder\'{o}n-Zygmund  decomposition to $f$ at height  $\alpha$ to obtain that there is a disjoint  family
of dyadic cubes $\{Q\}$ with  total  measure
$$\dsum_Q|Q| \lesssim{\alpha}^{-1}\|f\|_{L^1},$$
and that
$f=g+h,$ with
$\|g\|_{L^\infty}\lesssim\alpha,$
and $\|g\|_{L^1}\lesssim\|f\|_{L^1}$,
$h=\sum_{Q}h_Q,$ where
\text{supp}$(h_Q)\subseteq Q,$
$\int_{\mathbb R^n} h_Q(x)\,dx=0$ and
$\sum\|h_Q\|_{L^1}\lesssim\|f\|_{L^1}.$
By  Chebychev's inequality and \eqref{w3}, we get \begin{align*}\big| \big\{x: \Big(\sum_{k\in \Bbb Z}|({{\phi}}_k-\delta_0)\ast K^k\ast
g(x)|^q\Big)^{1/q}> \alpha\big\}\big|&\lesssim  \frac{\big\|\big(\sum_{k\in \Bbb Z}|({{\phi}}_k-\delta_0)\ast K^k\ast
g|^q\big)^{1/q}\big\|_{L^{2}}}{\alpha^2}\\&\lesssim  ( C_K+\|\omega\|_{Dini}+\|\omega\|_{Dini}^{1/q'}\|\omega^{1/2}\|_{Dini}^{2/q})\frac{\|g\|_{L^2}}{\alpha^2}\\&\lesssim (
C_K+\|\omega\|_{Dini}+\|\omega\|_{Dini}^{1/q'}\|\omega^{1/2}\|_{Dini}^{2/q})\frac{\|f\|_{L^1}}{\alpha}.\end{align*}
Then it suffices to estimate  $(\sum_{k\in \Bbb Z}|({{\phi}}_k-\delta_0)\ast K^k\ast h(x)|^q)^{1/q}$
away  from $\bigcup \widetilde{Q}$ with
 $\widetilde{Q}=2Q.$
Now by Chebychev's inequality and \eqref{jkx},
\begin{align*}&\alpha\big|\{x\notin\cup \widetilde{Q}:\Big(\sum_{k\in \Bbb Z}|({{\phi}}_{k}-\delta_0)\ast K^k\ast h(x)|^q\Big)^{1/q}>\alpha\}\big|\\
&\le\dsum_{Q}\dint_{x\notin \widetilde{Q}}\dint_Q|h_Q(y)|\Big(\sum_{k\in \Bbb Z}\big|(\phi_k-\delta_0)\ast K^k(x-y)-(\phi_k-\delta_0)\ast
K^k(x-y_Q)\big|^q\Big)^{1/q}dy\,dx\\
&\lesssim\dsum_{Q}\dint_Q|h_Q(y)|\dint_{|x-y_Q|\ge 2\ell(Q)}\bigg(\frac{C_K
|y-y_Q|^{\theta/2}}{|x-y_Q|^{n+\theta/2}}+\frac{\omega(|y-y_Q|/|x-y_Q|)^{\frac{1}{q'}}}{|x-y_Q|^{n}}\|\omega\|_{Dini}^{1/q}\\&\quad+\sum_{\frac{{2^k}}{2}\ge
|y-y_Q|}\frac{C_K}{|x-y_Q|^n}\chi_{2^{k}-|y-y_Q|\le|x-y_Q|\le2^{k}+|y-y_Q|}\bigg)\,dxdy\\
&\lesssim (C_K+\|{\omega}^{1/q'}\|_{Dini}\|\omega\|_{Dini}^{1/q})\dsum_{Q}\|h_Q\|_{L^1}\\&\quad+C_K\dsum_{Q}\dint_Q|h_Q(y)|\sum_{\frac{2^{k}}{2}\geq|y-y_Q|}
\frac{1}{2^k}\int_{2^{k}-|y-y_Q|\leq|x-y_Q|\leq2^{k}+|y-y_Q|}\frac{1}{|x-y_Q|^{n-1}}dx
dy\\
&\lesssim
(C_K+\|{\omega}^{1/q'}\|_{Dini}\|\omega\|_{Dini}^{1/q})\dsum_{Q}\|h_Q\|_{L^1}+C_K\dsum_{Q}\dint_Q|h_Q(y)
|\sum_{\frac{2^{k}}{2}\geq|y-y_Q|}\frac{|y-y_Q|}{2^{k}}dy\\
&\lesssim   (C_K+\|{\omega}^{1/q'}\|_{Dini}\|\omega\|_{Dini}^{1/q})\|f\|_{L^1},\end{align*}where $y_Q$ denotes the center of $Q.$ This finishes the proof of
\eqref{wj1}.

\quad

Now, we verify  \eqref{wj2}. Let $Q$ be a cube, and take $x,\,\xi\in Q.$ Let $B(x)=B(x,2\sqrt{n}\ell(Q)),$ then $3Q\subset B_x.$
By the triangle inequality, we get \begin{align*}
& \big(\sum_{k\in \Bbb Z}|(\phi_{k}-\delta_0)\ast K^k\ast (f{\chi_{{\Bbb R}^n\setminus 3Q}})(\xi)|^q\big)^{1/q}\\& \le \big|\big(\sum_{k\in \Bbb
Z}|(\phi_{k}-\delta_0)\ast K^k\ast (f{\chi_{{\Bbb R}^n\setminus B_x}})(\xi)|^q\big)^{1/q}-\big(\sum_{k\in \Bbb Z}|(\phi_{k}-\delta_0)\ast K^k\ast
(f{\chi_{{\Bbb R}^n\setminus B_x}})(x)|^q\big)^{1/q}\big|\\&\quad+\big(\sum_{k\in \Bbb Z}|(\phi_{k}-\delta_0)\ast K^k\ast f{\chi_{B_x\setminus
3Q}})(\xi)|^q\big)^{1/q}+\big(\sum_{k\in \Bbb Z}|(\phi_{k}-\delta_0)\ast K^k\ast (f{\chi_{{\Bbb R}^n\setminus B_x}})(x)|^q\big)^{1/q}\\&=:I+II+III.
\end{align*}
 For the term $I,$ by using the triangular inequality,
\begin{align}
I&\le\bigg( \sum_{k\in \Bbb Z}(\int_{{\Bbb R}^n \setminus
B_x}|(\phi_{k}-\delta_{0})*K^{k}(\xi-y)-(\phi_{k}-\delta_{0})*K^{k}(x-y)||f(y)|\,dy)^q\bigg)^{1/q}\\&\lesssim \bigg( \sum_{|x-\xi|\ge \frac{2^k}{2}}(\int_{{\Bbb R}^n
\setminus B_x}|(\phi_{k}-\delta_{0})*K^{k}(\xi-y)-(\phi_{k}-\delta_{0})*K^{k}(x-y)||f(y)|\,dy)^q\bigg)^{1/q}\nonumber\\&\quad+\bigg( \sum_{|x-\xi|\le
\frac{2^k}{2}}(\int_{{\Bbb R}^n \setminus B_x}|(\phi_{k}-\delta_{0})*K^{k}(\xi-y)-(\phi_{k}-\delta_{0})*K^{k}(x-y)||f(y)|\,dy)^q\bigg)^{1/q}\nonumber
\\&=:I_1+I_2.\nonumber
\end{align}
By using the size estimate \eqref{uj0} and the fact that $\frac{1}{2}|\xi-y|\le|x-y|\le \frac{3}{2}|\xi-y|$ (since $2|x-\xi|\le |x-y|$), we get
\begin{align}
I_1&\lesssim\bigg( \sum_{|x-\xi|\ge \frac{2^k}{2}}\big(\int_{{\Bbb R}^n \setminus B_x}\frac{\omega(2^k/|x-y|)}{|x-y|^{n}}\chi_{|x-y|>
\frac{3}{4}\cdot2^{k}}|f(y)|\,dy\big)^q\bigg)^{1/q}\\&\quad +C_K\int_{{\Bbb R}^n \setminus B_x}\sum_{|x-\xi|\ge \frac{2^k}{2}}\frac{|f(y)|}{|x-y|^n}\chi_{\frac{3}{4}\cdot2^{k}\le |x-y|\le
\frac{5}{4}\cdot2^{k}}\,dy\nonumber\\&\lesssim\bigg( \sum_{|x-\xi|\ge \frac{2^k}{2}}\big(\sum_{j\ge 1}\int_{2^{j}\sqrt{n}\ell(Q)\le|x-y|\le
2^{j+1}\sqrt{n}\ell(Q)}\omega\big(\frac{2^k}{2^{j}\sqrt{n}\ell(Q)}\big)\frac{1}{(2^{j}\sqrt{n}\ell(Q))^{n}}|f(y)|\,dy\big)^q\bigg)^{1/q}\nonumber\\&\quad+
C_K\int_{{\Bbb R}^n \setminus B_x}\frac{|x-\xi|}{|x-y|^{n+1}}|f(y)|\,dy\nonumber\\&\lesssim\Big( \sum_{|x-\xi|\ge
\frac{2^k}{2}}\omega\big(\frac{2^k}{\sqrt{n}\ell(Q)}\big)\Big)^{1/q}\sum_{j\ge 1}\omega^{1/q'}(2^{-j})Mf(x)+C_KMf(x)\nonumber\\&\lesssim
(\|\omega\|_{Dini}^{1/q}\|\omega^{1/q'}\|_{Dini}+C_K)Mf(x)\nonumber.
\end{align}
For the term $I_2$, by  \eqref{j1}, we get
\begin{align}
I_2&\lesssim\bigg( \sum_{|x-\xi|\le \frac{2^k}{2}}\big(\int_{{\Bbb R}^n \setminus B_x}\big(\frac{\omega(|x-\xi|/|x-y|)}{|x-y|^{n}}+
C_K\frac{|x-\xi|^\theta}{|x-y|^{n+\theta}}\big)\chi_{|x-y|\ge2^{k}/4}|f(y)|\,dy\big)^q\bigg)^{1/q}\nonumber\\&\quad+C_K\bigg( \sum_{|x-\xi|\le
\frac{2^k}{2}}\Big(\int_{{\Bbb R}^n \setminus
B_x}\frac{|f(y)|}{|x-y|^n}\chi_{2^{k}-|x-\xi|\leq|x-y|\le2^{k}+|x-\xi|}\,dy\Big)^q\bigg)^{1/q}\nonumber\\&\lesssim\bigg( \sum_{|x-\xi|\le \frac{2^k}{2}}\omega\big(\frac{|x-\xi|}{2^k}\big)+\big(\frac{|x-\xi|}{2^k}\big)^{\theta/2}\bigg)^{1/q}\int_{{\Bbb R}^n \setminus B_x}\bigg(\frac{\omega^{1/q'}(|x-\xi|/|x-y|)}{|x-y|^{n}}+
C_K\frac{|x-\xi|^{\theta/2}}{|x-y|^{n+\theta/2}}\bigg)|f(y)|\,dy\nonumber\\&\quad+\widetilde{I}_2\nonumber\\&\lesssim
(\|\omega\|_{Dini}^{1/q}\|\omega^{1/q'}\|_{Dini}+C_K)Mf(x)+\widetilde{I}_2\nonumber.
\end{align}
To estimate  $\widetilde{I}_2$, it suffices to consider the form
\begin{align*}
C_K\bigg(\sum_{|x-\xi|\le \frac{2^k}{2}}(\int_{s_{k}\le|x-y|\le s_{k+1}}\frac{1}{|x-y|^n}| f(y)|\,dy)^{q}\bigg)^{1/q},
\end{align*}
where
$|s_{k+1}-s_{k}| \leq 2|x-\xi|,$ and $\frac{2^{k}}{2}\le s_k\le \frac{3}{2}\cdot2^{k}.$ Using the hypothesis that $|x-\xi|<\ell(Q)$ and the kernel estimate we can bound the
above by a dimensional constant times
\begin{align*}
C_{K}\bigg(\sum_{|x-\xi|\le \frac{2^k}{2}}(s_{k}^{-n} \int_{s_{k}\le|x-y|\le s_{k+1}}|f(y)|\,dy)^{q}\bigg)^{1/q}
\end{align*}
The above $\ell^{q}$ norm can be written as
\begin{align*}
&~\bigg(\sum_{|x-\xi|\le \frac{2^k}{2}}(s_{k}^{-n}\int_{|x-y|\le s_{k+1}}|f(y)|\,dy-s_{k}^{-n}\int_{|x-y|\le s_{k}}|f(y)|\,dy)^{q}\bigg)^{1 / q} \\
&\leq\bigg(\sum_{k\in \Bbb Z}(s_{k+1}^{-n} \int_{|x-y|\le s_{k+1}}|f(y)|\,dy-s_{k}^{-n} \int_{|x-y|\le s_{k}}|f(y)|\,dy)^{q}\bigg)^{1/q}+\bigg(\sum_{|x-\xi|\le
\frac{2^k}{2}}((s_{k}^{-n}-s_{k+1}^{-n}) \int_{|x-y|\le s_{k+1}}|f(y)|\,dy)^{q}\bigg)^{1/q} \\
&\lesssim {V}_{q}\mathcal{A}(|f|)(x)+Mf(x)\bigg(\sum_{|x-\xi|\le \frac{2^k}{2}}((s_{k}^{-n}-s_{k+1}^{-n}) / s_{k+1}^{-n})^{q}\bigg)^{1/q},
\end{align*}
where $V_{q}\mathcal{A} (|f|)(x)=V_q\{A_{t}(|f|)(x)\}_{t>0}, $ and\begin{align*}
A_{t}(f)(x)=\frac{1}{|B_{t}|} \int_{B_{t}} f(x+y) dy, \quad x \in \mathbb{R}^{n}, t>0.
\end{align*}
Here $B_{t}$ denotes the open ball in $\mathbb{R}^{n}$ of center at the origin and radius $t$.
Also note that
\begin{align*}
\bigg(\sum_{|x-\xi|\le \frac{2^k}{2}}((s_{k}^{-n}-s_{k+1}^{-n}) / s_{k+1}^{-n})^{q}\bigg)^{1/q}
&\lesssim \sum_{|x-\xi|\le \frac{2^k}{2}} \frac{s_{k+1}-s_{k}}{2^k}
\lesssim \sum_{|x-\xi|\le \frac{2^k}{2}} \frac{|x-\xi|}{2^k}
 \le c_n.
\end{align*}
Thus   we get \begin{align}
\widetilde{I}_2&\lesssim
C_K({V}_{q}\mathcal{A}(|f|)(x)+Mf(x)).
\end{align}
Thus combining  the estimates of $I_1$ and $I_2,$ we get \begin{align*}
I&\le(\|\omega\|_{Dini}^{1/q}\|\omega^{1/q'}\|_{Dini}+C_K)Mf(x)+C_K {V}_{q}\mathcal{A}(|f|)(x).
\end{align*}

\quad

For the term $II$,     by using \eqref{tk1} and the fact that $|x-y|\simeq |\xi-y|$ (since $2|x-\xi|\le |x-y|$), we get
\begin{align}\label{I2}
II&\le\int_{{\Bbb R}^n} \sum_{k\in \Bbb Z}|(\phi_{k}-\delta_{0})*K^{k}(\xi-y)||f|{\chi_{B_x\setminus 3Q}}(y)\,dy\\&\lesssim
(\|\omega\|_{Dini}+C_K)\int_{3\ell(Q)\le|x-y|\le 2\sqrt{n}\ell(Q)}|\xi-y|^{-n}|f(y)|\,dy\nonumber \\&\lesssim(\|\omega\|_{Dini}+C_K)Mf(x).\nonumber
\end{align}

\quad

For the term $III$, since supp $(\phi_{k}-\delta_{0})*K^{k}\subset \{x:|x|\ge \frac{3}{4}\cdot2^k\},$ we obtain that
\begin{align*}
III
&\leq \bigg(\sum_{\frac{3}{4}\cdot2^k\geq2\sqrt{n} l(Q)}|(\phi_{k}-\delta_{0})*K^{k}*f\chi_{{\Bbb R}^n\setminus B_x}(x)|^{q}\bigg)^{1/q}+
\Bigg(\sum_{\frac{3}{4}\cdot2^k\leq2\sqrt{n} l(Q)}|(\phi_{j}-\delta_{0})*K^{k}*f\chi_{{\Bbb R}^n\setminus B_x}(x)|^{q}\Bigg)^{1/q}\\&= \bigg(\sum_{\frac{3}{4}\cdot2^k\geq2\sqrt{n}
l(Q)}|(\phi_{k}-\delta_{0})*K^{k}*f(x)|^{q}\bigg)^{1/q}+
\bigg(\sum_{\frac{3}{4}\cdot2^k\leq2\sqrt{n} l(Q)}|(\phi_{k}-\delta_{0})*K^{k}*f\chi_{{\Bbb R}^n\setminus B_x}(x)|^{q}\bigg)^{1/q}\\&\le \bigg(\sum_{k\in \Bbb
Z}|(\phi_{k}-\delta_{0})*K^{k}*f(x)|^{q}\Bigg)^{1/q}+
\Bigg(\sum_{\frac{3}{4}\cdot2^k\leq2\sqrt{n} l(Q)}|(\phi_{k}-\delta_{0})*K^{k}*f\chi_{{\Bbb R}^n\setminus B_x}(x)|^{q}\Bigg)^{1/q}.
\end{align*}
By using the size estimate \eqref{uj0} we get
\begin{align*}
|(\phi_{k}-\delta_{0})*K^{k}*f\chi_{\mathbb{R}^{n}\backslash B_{x}}(x)|
&\lesssim \int_{|x-y|>2\sqrt{n} l(Q)}\Big(\omega(2^k/|x-y|)+\frac{C_K2^k}{|x-y|}\chi_{\frac{3}{4}\cdot2^k\le |x-y|\le \frac{5}{4}\cdot2^k}\Big)\frac{|f(y)|}{|x-y|^n}dy\\
&\lesssim \sum_{j\ge 1}\int_{2^{j}\sqrt{n} l(Q)\le|x-y|\le2^{j+1}\sqrt{n}
l(Q)}\bigg(\omega(\frac{2^k}{2^j\ell(Q)})+C_K\frac{2^k}{2^j\ell(Q)}\bigg)\frac{|f(y)|}{(2^j\ell(Q))^n}dy\\&
\lesssim \Big(\sum_{j\ge 1}\omega(2^{-j})^{1/q'}\omega(\frac{2^k}{\ell(Q)})^{1/q}+C_K\frac{2^k}{\ell(Q)}\Big)Mf(x)\\&
\lesssim \bigg(\|\omega^{1/q'}\|_{Dini}\omega(\frac{2^k}{\ell(Q)})^{1/q}+C_K\frac{2^k}{\ell(Q)}\bigg)Mf(x),
\end{align*}
which implies
\begin{align*}
\bigg(\sum_{2^{k}\leq2\sqrt{n} l(Q)}|(\phi_{k}-\delta_{0})*K^{k}*f\chi_{\mathbb{R}^{n}\backslash B_{x}}(x)|^{q}\bigg)^{1/q}
&\lesssim  \bigg(\|\omega\|_{Dini}^{1/q}\|\omega^{1/q'}\|_{Dini}+C_K\bigg)Mf(x).
\end{align*}
Thus, we obtain that
\begin{align*}
III
&\lesssim  \Bigg(\sum_{k\in\mathbb{Z}}|(\phi_{k}-\delta_{0})*K^{k}*f(x)|^{q}\Bigg)^{1/q}+ (\|\omega\|_{Dini}^{1/q}\|\omega^{1/q'}\|_{Dini}+C_K)Mf(x).
\end{align*}

\quad

Combining  the estimates of $I,\,II$ and $III,$ we get
\begin{align*}&\mathcal{ M}_{(\sum_{j\in \Bbb Z}|(\phi_{k}-\delta_0)\ast K^k*f(x)|^{q})^{1/q}}\\&\lesssim  \bigg(\sum_{k\in \Bbb Z}|(\phi_{k}-\delta_0)\ast
K^k*f(x)|^{q}\bigg)^{1/q}+(C_K+ \|\omega\|_{Dini}^{1/q}\|\omega^{1/q'}\|_{Dini})Mf(x)+C_K {V}_{q}\mathcal{A}(|f|)(x).
\end{align*}
Then by the weak type $(1,1)$ of $M$ (see \cite{S1993}) and ${V}_{q}\mathcal{A}(|f|)$(see \cite{JRW03}), and \eqref{wj1}, we get
\begin{align*} \|\mathcal{M}_{(\sum_{k\in \Bbb Z}|(\phi_{k}-\delta_0)\ast K^k*f|^{q})^{1/q}}\|_{L^{1,\infty}}\lesssim  (C_K+
\|\omega\|_{Dini}^{1/q}\|\omega^{1/q'}\|_{Dini}+\|\omega\|_{Dini}^{1/q'}\|\omega^{1/2}\|_{Dini}^{2/q}+\|\omega\|_{Dini})\|f\|_{L^{1}},
\end{align*} which gives the proof of  \eqref{wj2}.

The proof of Proposition \ref{pro4} is complete.\qed

\quad

\section{ Proofs of Theorem \ref{thm3} and Theorem \ref{thm1}}
\emph{Proof of Theorem \ref{thm1}.}    Recalling that for any  fixed $k\in \Bbb Z,$ we denote by
\begin{align*}K_{k}(x)=K(x)\chi_{|x|\leq2^{k}}\end{align*}
and\begin{align*}K^{k}(x)=K(x)\chi_{|x|>2^{k}}.\end{align*}
Then for any fixed $k\in \Bbb Z,$ write
\begin{align*}
\phi_{k}*K(x)=\phi_{k}*K_{k}(x)+(\phi_{k}-\delta_{0})*K^{k}(x)+K^{k}(x),
\end{align*}where $\delta_0$ is the Dirac measure at $0.$
Thus by the triangle inequality, we get
\begin{align*}
&V_{q}\{\phi_{k}*K*f\}_{k\in \mathbb{Z}}\\&\leq V_{q}\{\phi_{k}*K_{k}*f\}_{k\in \mathbb{Z}}+V_{q}\{(\phi_{k}-\delta_{0})*K^{k}*f\}_{k\in
\mathbb{Z}}+V_{q}\{K^{k}*f\}_{k\in \mathbb{Z}}\\
&\leq \Bigg(\sum_{k\in \mathbb{Z}}|\phi_{k}*K_{k}*f|^{q}\Bigg)^{1/q}+\Bigg(\sum_{k\in \mathbb{Z}}|(\phi_{k}-\delta_{0})*K^{k}*f|^{q}\Bigg)^{1/q}+
V_{q}\{K^{k}*f\}_{k\in \mathbb{Z}}.
\end{align*}
Since the following results have been established in \cite {deZ16, DLY19, CM},
\begin{align}\label{vk2}
\|V_{q}\{K^{k}*f\}_{k\in \mathbb{Z}}\|_{L^{1,\infty}}&\lesssim (\|\omega\|_{Dini}+C_K+\|T\|_{L^2\rightarrow L^2})\|f\|_{L^{1}}.
\end{align} and for $1<p<\infty,$ $w\in A_p,$
\begin{align}\label{vk3}
\|V_{q}\{K^{k}*f\}_{k\in \mathbb{Z}}\|_{L^{p}(w)}&\lesssim (\|\omega\|_{Dini}+C_K+\|T\|_{L^2\rightarrow L^2})\{w\}_{A_p}\|f\|_{L^{p}(w)}.
\end{align}
Then combining Proposition \ref{pro2} and Proposition \ref{pro4}, we finish the proof of Theorem \ref{thm1}. \qed

\quad\quad

\emph{Proof of Theorem \ref{thm3}.}     Recall the definition of the operator $T_{\Omega}$ given in the introduction. It can be
written as
\begin{align}\label{eq-01}
T_{\Omega}f=\sum_{k \in \mathbb{Z}} T_{k} f=\sum_{k \in \mathbb{Z}} v_{k}\ast f, \quad v_{k}=\frac{\Omega(x^{\prime})}{|x|^{n}}
\chi_{\{2^{k}<|x|\le2^{k+1}\}}.
\end{align}

We consider the following partition of unity. Let $\varphi \in C_{c}^{\infty}(\mathbb{R}^{n})$ be such that
 supp $\varphi \subset\{x:|x| \leq \frac{1}{100}\}$ and $\int \varphi d x=1,$ and so that $\widehat{\varphi} \in \mathcal{S}(\mathbb{R}^{n}) .$ Let us also
define
$\psi$ by $\widehat{\psi}(\xi)=\widehat{\varphi}(\xi)-\widehat{\varphi}(2 \xi) .$ Then, with this choice of $\psi,$ it follows that $\int \psi d x=0$.
We write $\varphi_{j}(x)=\frac{1}{2^{j n}} \varphi(\frac{x}{2^{j}}),$ and $\psi_{j}(x)=\frac{1}{2^{j n}} \psi(\frac{x}{2^{j}})$.
We now define the partial sum operators $S_{j}$ by $S_{j}(f)=f\ast\varphi_{j}$. Their differences
are given by
\begin{align}\label{eq-02}
S_{j}(f)-S_{j+1}(f)=f\ast\psi_{j}.
\end{align}
Since $S_{j}f \rightarrow f$ as $j \rightarrow -\infty$, for any sequence of integers $\mathcal{N}=\{N(j)\}_{j=0}^{\infty}$, with
$0 = N(0) < N(1) <\cdot\cdot\cdot < N(j)\rightarrow\infty $, we have the identity
\begin{equation}\label{eq-03}
T_{k}=T_{k} S_{k}+\sum_{j=1}^{\infty} T_{k}(S_{k-N(j)}-S_{k-N(j-1)}).
\end{equation}
In this way, \begin{equation}\label{eq-03}T_{\Omega}=\sum_{j=0}^{\infty} \widetilde{T}_{j}=\sum_{j=0}^{\infty} {T}_{j}^{\mathcal{N}},\end{equation} where
\begin{equation}\label{eq-04}
\widetilde{T}_{0}:={T}_{0}^{\mathcal{N}}:=\sum_{k \in \mathbb{Z}} T_{k} S_{k}
\end{equation}
and, for $j\geq1$,

\begin{align}\label{eq-05} \widetilde{T}_{j} &:=\sum_{k \in \mathbb{Z}} T_{k}(S_{k-j}-S_{k-(j-1)}),\nonumber\\
{T}_{j}^{\mathcal{N}}
&:=\sum_{k \in \mathbb{Z}} T_{k}(S_{k-N(j)}-S_{k-N(j-1)})=\sum_{i=N(j-1)+1}^{N(j)} {\widetilde{T}}_{i}. \end{align}
Then
\begin{align*}V_q\{\phi_l\ast T_\Omega f\}_{l\in \Bbb Z}&=V_q\big\{\phi_l\ast (\sum_{j=0}^{\infty}{T}_{j}^{\mathcal{N}} f)\big\}_{l\in \Bbb Z}
\leq\sum_{j=0}^{\infty} V_q\{\phi_l\ast{T}_{j}^{\mathcal{N}}f\}_{l\in \Bbb Z}. \end{align*}
Therefore, \begin{align*}\|V_q\{\phi_l\ast T_\Omega f\}_{l\in \Bbb Z}\|_{L^p(w)}\leq\sum_{j=0}^{\infty} \|V_q\{\phi_l\ast{T}_{j}^{\mathcal{N}}f\}_{l\in
\Bbb Z}\|_{L^p(w)}. \end{align*}
To prove Theorem \ref{thm3}, we claim that the following inequalities hold for $1<p<\infty$ and $w\in A_p$:
\begin{align}\label{vv}
\|V_q\{\phi_{l}\ast {T}_{j}^{\mathcal{N}}f\}_{l\in \Bbb Z}\|_{L^p}&\lesssim \|\Omega\|_{L^{\infty}}2^{-\theta
N(j-1)}(1+N(j))^{1+1/q}\|f\|_{L^p}
 \end{align}
 and
\begin{align}\label{vq3}
\|V_q\{\phi_{l}\ast {T}_{j}^{\mathcal{N}}f\}_{l\in \Bbb Z}\|_{L^p(w)}&\lesssim
\|\Omega\|_{L^{\infty}}(1+N(j))^{1+1/q}\{w\}_{A_p}\|f\|_{L^p(w)}.
\end{align}

Assume the above two claims at the moment.
Set $\varepsilon:=\frac{1}{2}c_n/(w)_{A_p},$ $c_n$ is small enough (see \cite[Corollary 3.18]{HRT}). By \eqref{vq3}, we have, for this choice of $\varepsilon$,
\begin{align}\label{vw}
\|V_q(\{\phi_{l}\ast {T}_{j}^{\mathcal{N}}f\}_{l\in \Bbb
Z})\|_{L^p(w^{1+\varepsilon})}&\lesssim\|\Omega\|_{L^\infty}(1+N(j))^{1+1/q}\{w^{1+\varepsilon}\}_{A_p}\|f\|_{L^p(w^{1+\varepsilon})}\\
&\lesssim\|\Omega\|_{L^\infty}(1+N(j))^{1+1/q}\{w\}_{A_p}^{1+\varepsilon}\|f\|_{L^p(w^{1+\varepsilon})}\nonumber.
\end{align}

Now we are in position to apply the interpolation theorem with change of measures by Stein and Weiss---Lemma \ref{sw}.
We apply it to $T=V_q\{\phi_{l}\ast {T}_{j}^{\mathcal{N}}\}_{l\in \Bbb Z}$ with $p_0=p_1=p,~~w_0=w^0=1$ and $w_1=w^{1+\varepsilon}$, so that by
$\lambda=\varepsilon/(1+\varepsilon)$, \eqref{vv} and \eqref{vw}, one has for some $\theta,\,\gamma>0 $ such that
\begin{align*}
\|V_q\{\phi_{l}\ast {T}_{j}^{\mathcal{N}}\}_{l\in \Bbb Z}\|_{L^p(w)\rightarrow L^p(w)}&\lesssim
\|V_q\{\phi_{l}\ast {T}_{j}^{\mathcal{N}}\}_{l\in \Bbb Z}\|_{L^p\rightarrow L^p}^{\varepsilon/(1+\varepsilon)}\|V_q\{\phi_{l}\ast
{T}_{j}^{\mathcal{N}}\}_{l\in \Bbb Z}\|_{L^p(w^{1+\varepsilon})\rightarrow L^p(w^{1+\varepsilon})}^{1/(1+\varepsilon)}\\
&\lesssim\|\Omega\|_{L^\infty}(1+N(j))^{1+1/q}2^{-\theta N(j-1)\varepsilon/(1+\varepsilon)}\{w\}_{A_p}\\
&\lesssim\|\Omega\|_{L^\infty}(1+N(j))^{1+1/q}2^{-\gamma N(j-1)/{(w)}_{A_p}}\{w\}_{A_p}.
\end{align*}
Thus
\begin{align*}
\|V_q\{\phi_{l}\ast T_\Omega\}_{l\in \Bbb Z}\|_{L^p(w)\rightarrow L^p(w)}&\lesssim
\sum^\infty_{j=0}\|V_q\{\phi_{l}\ast {T}_{j}^{\mathcal{N}}\}_{l\in \Bbb Z}\|_{L^p(w)\rightarrow L^p(w)}\\
&\lesssim \|\Omega\|_{L^\infty}\{w\}_{A_p}\sum^\infty_{j=0}(1+N(j))^{1+1/q}2^{-\gamma N(j-1)/{(w)}_{A_p}}.
\end{align*}
 We now choose $N(j)=2^j$ for $j\geq1$. Then, using $e^{-x}\leq2x^{-2}$, we
have
\begin{align*}
\sum^\infty_{j=0}&(1+N(j))^{1+1/q}2^{-\gamma N(j-1)/{(w)}_{A_p}}
\lesssim\sum_{j:2^{j}\leq(w)_{A_p}}2^{j(1+1/q)}+\sum_{j:2^j\geq(w)_{A_p}}2^{j(1+1/q)}\big(\frac{(w)_{A_p}}{2^j}\big)^2
\lesssim(w)_{A_p}^{1+1/q},
\end{align*}
by summing the two geometric series in the last step.
This implies
$$
\|V_q\{\phi_{l}\ast T_\Omega f\}_{l\in \Bbb Z}\|_{L^p(w)}\lesssim\|\Omega\|_{L^\infty}\{w\}_{A_p}(w)_{A_p}^{1+1/q}\|f\|_{L^p(w)},
$$
and hence the proof of Theorem \ref{thm3} is complete under the assumptions that \eqref{vv} and \eqref{vq3} hold.
\quad

Now we turn to proving the claims \eqref{vv} and \eqref{vq3}.
The following inequality  can be found in \cite[Lemma 3.7]{HRT},
\begin{align}\label{jn}
\|{T}_{j}^{\mathcal{N}}f\|_{L^{2}}
& \lesssim \|\Omega\|_{L^{\infty}} 2^{-\alpha N(j-1)}\|f\|_{L^{2}},
\end{align}
for some  $0 < \alpha < 1$ independent of $T_{\Omega}$ and $j.$
Then by  (see \cite{JSW08, DHL})\begin{align}
\|V_q\{\phi_l\ast f\}_{l\in \Bbb Z}\|_{L^2}\lesssim\|f\|_{L^2},
\end{align}
we get
\begin{align}\label{v1}\|V_q\{\phi_l\ast{T}_{j}^{\mathcal{N}}f\}_{l\in \Bbb Z}\|_{L^2} \lesssim  \|{T}_{j}^{\mathcal{N}}f\|_{L^{2}}\lesssim
\|\Omega\|_{L^{\infty}} 2^{-\alpha N(j-1)}\|f\|_{L^{2}}. \end{align}

The operator ${T}_{j}^{\mathcal{N}} $ is a $\omega$-Calder\'{o}n--Zygmund operator with the kernel $K_j^\mathcal{N}$ (see \cite[Lemma 3.10]{HRT}) satisfying
\begin{align} \label{Kj1}|K_j^\mathcal{N}(x)|  \lesssim \frac{\|\Omega\|_{L^{\infty}}}{|x|^n}  \end{align}and for $2|y|\le |x|,$ \begin{align}\label{Kj2}
|K_j^\mathcal{N}(x-y)-K_j^\mathcal{N}(x)|  \lesssim \frac{\omega_j(|y|/|x|)}{|x|^n}, \end{align} where   $\omega_{j}(t)\le\|\Omega\|_{L^ \infty} \min
(1,2^{N(j)}t), $ and $\|\omega_j\|_{Dini}\lesssim\|\Omega\|_{L^{\infty}}(1+N(j)). $ For $r>1,$ the Dini norm of $\omega_{j}^{1/r}$ is estimated as  \begin{align}\|\omega_j^{1/r}\|_{Dini}\lesssim \|\Omega\|_{L^{\infty}}^{1/r} \int_0^{2^{-N(j)}}2^{N(j)/r}t^{1/r}\frac{dt}{t}+\|\Omega\|_{L^{\infty}}^{1/r}\int_{2^{-N(j)}}^{1}\frac{dt}{t}\lesssim\|\Omega\|_{L^{\infty}}^{1/r}(1+N(j)). \end{align}
Applying Theorem \ref{thm1} to $T={T}_{j}^{\mathcal{N}}$, we get
for $1<p<\infty$ and $w\in A_p$\begin{align*} \|V_q\{\phi_{l}\ast {T}_{j}^{\mathcal{N}}f\}_{l\in \Bbb Z}\|_{L^p(w)}&\lesssim
\|\Omega\|_{L^{\infty}}(1+N(j))^{1+1/q}\{w\}_{A_p}\|f\|_{L^p(w)}, \end{align*} which gives the proof of  \eqref{vq3}.
Taking $w=1,$ we get for $1<p<\infty,$ \begin{align}\label{tj2} \|V_q\{\phi_{l}\ast {T}_{j}^{\mathcal{N}}f\}_{l\in \Bbb Z}\|_{L^p}&\lesssim
\|\Omega\|_{L^{\infty}}(1+N(j))^{1+1/q}\|f\|_{L^p}. \end{align}
Interpolating between \eqref{v1} and \eqref{tj2}, we get\begin{align*} \|V_q\{\phi_{l}\ast {T}_{j}^{\mathcal{N}}f\}_{l\in \Bbb Z}\|_{L^p}&\lesssim
\|\Omega\|_{L^{\infty}}2^{-\theta N(j-1)}(1+N(j))^{1+1/q}\|f\|_{L^p} \end{align*} which establishes the proof of \eqref{vv}.

 The proof of Theorem \ref{thm3} is complete. \qed

\quad

\section{Proof of Theorem \ref{pro01}} Recalling that for any fixed $k\in \Bbb Z,$  $$T_{\Omega,2^k}f(x)=\int_{|x-y|>2^k}\frac{\Omega(x-y)}{|x-y|^n}f(y)\,dy,$$  and $K_\Omega(x)=\frac{\Omega(x')}{|x|^n}$. Following \cite{DR86}, for any fixed $k\in \Bbb Z,$ we write
\begin{align*}
T_{\Omega, 2^k}f(x)
&=\phi_{k}\ast T_\Omega f-\phi_{k}\ast K_\Omega\chi_{|\cdot|\le 2^k}\ast f+(\delta_0-\phi_{k})\ast T_{\Omega,2^k} f,
\end{align*}   which splits $\mathcal{T}_\Omega =\{T_{\Omega, 2^k}\}_{k\in \Bbb Z}$ into three families:
\begin{align*}
\mathcal{T}_\Omega^1(f):=\{\phi_{k}\ast T_\Omega f\}_{k\in \Bbb Z},\quad
\mathcal{T}_\Omega^2(f):=\{\phi_{k}\ast (K_\Omega\chi_{|\cdot|\leq 2^k})\ast f\}_{k\in \Bbb Z},\quad
\mathcal{T}_\Omega^3(f):=\{(\delta_0-\phi_{k})\ast T_{\Omega,2^k}f\}_{k\in \Bbb Z}.
\end{align*}
Thus, it suffices to estimate the weighted $L^p$ norm of $\mathcal{T}_\Omega^i(f)$, $i=1,2,3$.


\quad

\textbf{Part 1.} Let us first consider $\mathcal{T}_\Omega^1(f)$. By Theorem \ref{thm3}, we get that for $1<p<\infty$ and $w\in A_p,$
\begin{align*}
\|V_q(\mathcal{T}_\Omega^1(f))\|_{L^p(w)}\lesssim \|\Omega\|_{L^{\infty}} (w)_{A_p}^{1+1/q} \{w\}_{A_p}\|f\|_{L^p(w)}.
\end{align*}

\textbf{Part 2.}  For the term $\mathcal{T}_\Omega^2(f)$, by the Minkowski inequality,  we get for $1<p<\infty$ and $w\in A_p,$
\begin{align*}
\|V_q(\mathcal{T}_\Omega^2(f))\|_{L^p(w)}
\lesssim\bigg\|\bigg(\sum_{k\in \Bbb Z}|\phi_k\ast K_\Omega\chi_{|\cdot|\le 2^k}\ast f|^q\bigg)^{1/q}\bigg\|_{L^p(w)}.
\end{align*}
 It is easy to verify that $K_\Omega$ satisfies the cancellation condition \eqref{2} and the size estimate $|K_\Omega(x)|\le
\frac{\|\Omega\|_{L^\infty}}{|x|^n}.$ Then applying Proposition \ref{pro2} to $K=K_\Omega$, we get  \begin{align*}\bigg\|\bigg(\sum_{k\in \Bbb Z}|\phi_k\ast K_\Omega\chi_{|\cdot|\le 2^k}\ast
f|^q\bigg)^{1/q}\bigg\|_{L^p(w)}\lesssim \|\Omega\|_{L^\infty}\{w\}_{A_p}\|f\|_{L^p(w)},\end{align*} by $\ell^2\subset\ell^q.$

\textbf{Part 3.} We now estimate $\mathcal{T}_\Omega^3(f).$
We claim that  for $1<p<\infty$ and $w\in A_p,$
\begin{align}\label{pro03}
\|V_q(\mathcal{T}_\Omega^3(f))\|_{L^p(w)}\lesssim\|\Omega\|_{L^\infty}(w)_{A_p}^{1+1/q}\{w\}_{A_p}\|f\|_{L^p(w)},
\end{align}
which, together with the estimates of $\mathcal{T}_\Omega^1(f)$, $\mathcal{T}_\Omega^2(f)$ and $\mathcal{T}_\Omega^3(f),$ gives
\begin{align*}
&\|V_q\mathcal{T}_\Omega(f) \|_{L^p(w)}\lesssim\|\Omega\|_{L^\infty}(w)_{A_p}^{1+1/q}\{w\}_{A_p}\|f\|_{L^p(w)}.
\end{align*} This finishes the proof of Theorem \ref{pro01}.

\quad

We now provide the proof of \eqref{pro03}.
For $k\in \Bbb Z$, we define $v_{k}$ as
$$v_{k}(x)=\frac{\Omega(x')}{|x|^n}\chi_{2^k<|x|\le 2^{k+1}}(x).$$ Then we can write $T_{\Omega, 2^k}f(x)=\sum_{s\ge 0}v_{k+s}\ast f(x).$
Thus \begin{align*}\|V_q(\mathcal{T}_\Omega^3(f))\|_{L^p(w)}=\bigg\|V_q\bigg\{(\delta-\phi_{k})\ast\sum_{s\ge 0} v_{k+s}\ast
f\bigg\}\bigg\|_{L^p(w)}. \end{align*} Recall that $S_{j}(f)=f\ast\varphi_{j}$ and  their differences
are given by
\begin{align}\label{eq-0s}
S_{j}(f)-S_{j+1}(f)=f\ast\psi_{j},
\end{align} where $\varphi_{j}$ and $\psi_{j}$ are
defined as  in  Section 5.
Denote by $T_{k}f(x)=v_k\ast f(x).$ Similarly,  for any sequence of integers $\mathcal{N}=\{N(j)\}_{j=0}^{\infty}$, with
$0 = N(0) < N(1) <\cdot\cdot\cdot < N(j)\rightarrow\infty $, we also have the identity.
\begin{align}\label{ksj}
T_{k+s}&=T_{k+s}S_{k+s}+\sum_{j=1}^\infty T_{k+s}(S_{k+s-N(j)}-S_{k+s-N(j-1)})
=:\sum_{j=0}^\infty T^\mathcal{N}_{k+s,j},
\end{align}and\begin{align}\label{ksj'}
v_{k+s}
=:\sum_{j=0}^\infty K^\mathcal{N}_{k+s,j},
\end{align}
where $T_{k+s,0}^\mathcal{N}=T_{k+s}S_{k+s}$ and for $j\geq1$,
\begin{align}\label{kej}
 T^\mathcal{N}_{k+s,j}&=T_{k+s}(S_{k+s-N(j)}-S_{k+s-N(j-1)})\\ \nonumber&=\sum_{i=N(j-1)+1}^{N(j)}T_{k+s}(S_{k+s-i}-S_{k+s-(i-1)})\\
 \nonumber&=:\sum_{i=N(j-1)+1}^{N(j)}T_{k+s,i}.
\end{align}
Then
\begin{align}\label{vej}
(\delta_0-\phi_k)\ast v_{k+s}\ast f(x)
=(\delta_0-\phi_k)\ast T_{k+s} f(x)=\sum^\infty_{j=0}(\delta_0-\phi_k)\ast T^\mathcal{N}_{k+s,j}f(x).
\end{align}
So
\begin{align*}
\|V_q(\mathcal{T}_\Omega^3(f))\|_{L^p(w)}
&\leq\sum^\infty_{j=0}\bigg\|V_q\bigg(\bigg\{(\delta_0-\phi_k)\ast\sum_{s\geq0} T^\mathcal{N}_{k+s,j}f\bigg\}_{k\in \Bbb
Z}\bigg)\bigg\|_{L^p(w)}\\&\leq\sum^\infty_{j=0}\bigg\|\bigg(\sum_{k\in \Bbb Z}|(\delta_0-\phi_k)\ast\sum_{s\geq0}
T^\mathcal{N}_{k+s,j}f|^q\bigg)^{1/q}\bigg\|_{L^p(w)}.
\end{align*}
We first estimate the $L^2-$norm of $(\sum_{k\in \Bbb Z}|(\delta_0-\phi_k)\ast\sum_{s\geq0} T^\mathcal{N}_{k+s,j}f|^q)^{1/q}$.
\begin{lemma}  \label{ls1}   We have for $j\ge 0$, there is a positive $0<\tau<1$ such that
\begin{align}\label{phi3}\bigg\|\bigg(\sum_{k\in \Bbb Z}|(\delta_0-\phi_k)\ast\sum_{s\geq0} T^\mathcal{N}_{k+s,j}f|^q\bigg)^{1/q}\bigg\|_{ L^2}\lesssim
\|\Omega\|_{L^\infty}2^{-\tau N(j-1)} \|f\|_{L^2}.\end{align}
\end{lemma}
{\emph{Proof}.} By the Plancherel Theorem and  \eqref{kej},  we get for $j\ge 1,$ \begin{align*}
&\bigg\|\bigg(\sum_{k\in \Bbb Z}|(\delta_0-\phi_k)\ast \sum_{s\ge 0}T^\mathcal{N}_{k+s,j}f|^q\bigg)^{1/q}\bigg\|_{ L^2}\\&\le \sum_{s\geq0}\bigg\|\bigg(\sum_{k\in
\mathbb{Z}}|(\delta_0-\phi_k)\ast T^\mathcal{N}_{k+s,j}f|^2\bigg)^\frac{1}{2}\bigg\|_{L^2}\\
&=\sum_{s\geq0}\bigg(\sum_{k\in \mathbb{Z}}\int_{{\Bbb R}^n}|\widehat{(\delta_0-\phi_k)}(\xi)|^2|\widehat{T^\mathcal{N}_{k+s,j}f}(\xi)|^2d\xi\bigg)^{1/2}\\
&\le\sum_{s\geq0}\sum_{i=N(j-1)+1}^{N(j)}\bigg(\sum_{k\in \mathbb{Z}}\int_{{\Bbb
R}^n}|\widehat{(\delta_0-\phi_k)}(\xi)|^2|\widehat{T_{k+s,i}f}(\xi)|^2d\xi\bigg)^{1/2}.
\end{align*}
Since $\widehat{\phi}\in \mathcal{S}({\Bbb R}^n)$ and $\widehat{\phi}(0)=1,$ $|1-\widehat{\phi}(\xi)|\lesssim \min(1, |\xi|).$  Also by
$\widehat{\psi}\in \mathcal{S}({\Bbb R^n})$ and $\widehat{\psi}(0)=0,$ then $|\widehat{\psi}(\xi)|\lesssim \min(1, |\xi|).$ Therefore,
$|1-\widehat{\phi_k}(\xi)|\lesssim \min(1,|2^k\xi|^\alpha)$, $|\widehat{\psi}(2^{k+s-i}\xi)|\lesssim
\min(|2^{k+s-i}\xi|^\gamma,1),\,|\widehat{v}_{k+s}(\xi)|\lesssim \|\Omega\|_{L^\infty}|2^{k+s}\xi|^{-\beta}$ (see \cite{DR86}) for some
$0<\beta,\,\gamma,\,\alpha<1$, $\alpha+\gamma-\beta>0,$ $\gamma<\beta$ and $\alpha<\beta,$
\begin{align*}
&\sum_{k\in \mathbb{Z}}|\widehat{(\delta_0-\phi_k)}(\xi)|^2|\widehat{T_{k+s,i}f}(\xi)|^2\\&=\sum_{k\in
\mathbb{Z}}|1-\widehat{\phi_k}(\xi)|^2|\widehat{v}_{k+s}(\xi)|^2|\widehat{\psi}(2^{k+s-i}\xi)|^2|\widehat{f}(\xi)|^2\\
&\lesssim \|\Omega\|_{L^\infty}^2\bigg(\sum_{|2^k\xi|\leq2^i}|2^k\xi|^{2\alpha}|2^{k+s-i}\xi|^{2\gamma
}|2^{k+s}\xi|^{-2\beta}+\sum_{|2^k\xi|\geq2^i}|2^{k+s}\xi|^{-2\beta}\bigg)|\widehat{f}(\xi)|^2\\
&\lesssim \|\Omega\|_{L^\infty}^2\bigg(2^{-2\gamma i}2^{-2(\beta-\gamma)
s}\sum_{|2^k\xi|\leq2^i}|2^k\xi|^{2\alpha+2\gamma-2\beta}+2^{-2s\beta}\sum_{|2^k\xi|\geq2^i}|2^{k}\xi|^{-2\beta}\bigg)|\widehat{f}(\xi)|^2\\
&\lesssim\|\Omega\|_{L^\infty}^2\big(2^{-2(\beta-\alpha)i}2^{-2(\beta-\gamma) s}+2^{-s\beta}2^{-2\beta i}\big)|\widehat{f}(\xi)|^2\\
&\lesssim\|\Omega\|_{L^\infty}^22^{-2(\beta-\alpha)i}2^{-2(\beta-\gamma)s}|\widehat{f}(\xi)|^2.
\end{align*}
Therefore, by summing the geometric series and the Plancherel theorem, we get for some $\tau\in (0,1),$
\begin{align*}
&\bigg\|\bigg(\sum_{k\in \Bbb Z}|(\delta_0-\phi_k)\ast \sum_{s\geq0}T^\mathcal{N}_{k+s,j}f|^q\bigg)^{1/q}\bigg\|_{L^2}\lesssim\|\Omega\|_{L^\infty}2^{-\tau N(j-1)}
\|f\|_{L^2}.
\end{align*}
For $j=0$, $|1-\widehat{\phi}(2^{k}\xi)|\lesssim \min(|2^{k}\xi|,1),\,|\widehat{\varphi}(2^{k}\xi)|\lesssim1,\,|\widehat{v}_{k+s}(\xi)|\lesssim\|\Omega\|_{L^\infty}|2^{k+s}\xi|^{-\beta}$ for some
$0<\beta<1$ (see \cite{DR86}),
\begin{align*}
\sum_{k\in \mathbb{Z}}|\widehat{(\delta_0-\phi_k)}(\xi)|^2|\widehat{T_{k+s,0}f}(\xi)|^2&=\sum_{k\in
\mathbb{Z}}|1-\widehat{\phi_k}(\xi)|^2|\widehat{v}_{k+s}(\xi)|^2|\widehat{\varphi}(2^{k}\xi)|^2|\widehat{f}(\xi)|^2\\
&\lesssim\|\Omega\|_{L^\infty}^2\bigg(\sum_{|2^k\xi|\leq 1}|2^{k}\xi|^2|2^{k+s}\xi|^{-2\beta}+\sum_{|2^k\xi|\geq
1}|2^{k+s}\xi|^{-2\beta}\bigg)|\widehat{f}(\xi)|^2\\
&\lesssim\|\Omega\|_{L^\infty}^2\bigg(2^{-2\beta
s}\sum_{|2^k\xi|\leq1}|2^k\xi|^{2-2\beta}+2^{-2s\beta}\sum_{|2^k\xi|\geq1}|2^{k}\xi|^{-2\beta}\bigg)|\widehat{f}(\xi)|^2\\
&\lesssim\|\Omega\|_{L^\infty}^22^{-2\beta s}|\widehat{f}(\xi)|^2.
\end{align*}
Therefore, by the Plancherel theorem, we get for some $\tau\in (0,1),$
\begin{align*}
&\bigg\|\bigg(\sum_{k\in \Bbb Z}|(\delta_0-\phi_k)\ast\sum_{s\geq0} T^\mathcal{N}_{k+s,0}f|^2\bigg)^{1/2}\bigg\|_{L^2}\lesssim\|\Omega\|_{L^\infty} \|f\|_{L^2}.
\end{align*}

The proof of Lemma \ref{ls1} is complete by $\ell^2\subset \ell^q$. \qed

\bigskip
Denote $T_{k,j}^\mathcal{N}:=\sum_{s\geq0} T^\mathcal{N}_{k+s,j}$, and let
 $K_{k,j}^\mathcal{N}=\sum_{s\geq0} K^\mathcal{N}_{k+s,j}$ be the kernel of $T_{k,j}^\mathcal{N}$
given by \begin{align}\label{gs1}K_{k,j}^\mathcal{N}:=\sum_{s\geq0}
v_{k+s}\ast(\varphi_{{k+s-N(j)}}-\varphi_{{k+s-N(j-1)}})\end{align} for $j\ge 1$ and for $j=0,$
\begin{align}\label{gs0}K_{k,0}^\mathcal{N}:=\sum_{s\geq0}  v_{k+s}\ast \varphi_{{k+s}}.\end{align} First, it is easy to verify that
supp $K_{k,j}^\mathcal{N}\subset\{x: |x|\ge \frac{3}{4}\cdot2^{k}\}.$ In the following, we will verify  that  $K_{k,j}^\mathcal{N}$ is a $\omega$-Dini Calder\'{o}n-Zygmund kernel satisfying
\eqref{K0} and \eqref{K1}.

\quad

\begin{lemma}  \label{l1}For $j\ge 0$ and $k\in \Bbb Z$. Then we have the size estimate
\begin{align}\label{g8}
|K_{k,j}^\mathcal{N}(x)|
&\lesssim\frac{\|\Omega\|_{L^\infty}}{|x|^n}\chi_{
|x|\ge\frac{3}{4}\cdot2^{k}}(x).
\end{align}
and the regularity estimate
\begin{align}\label{gr8}
|K_{k,j}^\mathcal{N}(x-h)-K_{k,j}^\mathcal{N}(x)|
&\lesssim\frac{\omega_j(\frac{|h|}{|x|})}{|x|^n}\chi_{
|x|\ge\frac{3}{4}\cdot2^{k}}(x),\quad  0<|h|<\frac{|x|}{2},
\end{align}
where $\omega_j(t)\le\|\Omega_j\|_{L^{\infty}}\min(1, 2^{N(j)}t).$
\end{lemma}

 \quad

\emph{Proof.}  In order to get the required estimates for the kernel $K_{k,j}^\mathcal{N}$, we first study the kernel of each $H_{k, s, N(j)}$ which is
defined by $$H_{k, s, N(j)}:= v_{k+s}\ast \varphi_{{k+s-N(j)}}.$$ First, we estimate $|H_{k, s, N(j)}(x)|.$ Since supp $\varphi\subset \{x: |x|\le
\frac{1}{100}\}$, a simple computation gives that
\begin{align}\label{v5}
|H_{k,s,N(j)}(x)|
&\le\int_{2^{k+s}\leq |y|\leq2^{k+s+1}}\frac{|\Omega(y)|}{|y|^{n}}|\varphi_{{k+s-N(j)}}(x-y)|dy\\
&\lesssim\frac{\|\Omega\|_{L^{\infty}}}{|x|^{n}}\chi_{\frac{3}{4}\cdot2^{k+s}\leq |x|\leq\frac{5}{4}\cdot2^{k+s}}.\nonumber
\end{align}
 From the triangular inequality,  and $N(j-1)<N(j)$, we obtain that the kernel $K_{k,j}^\mathcal{N}:=\sum_{s\geq0}(H_{k,s,N(j)}-H_{k,s,N(j-1)})$ satisfies
 \begin{align}\label{h1}
|K_{k,j}^\mathcal{N}(x)|
&\lesssim \sum_{s\geq0}\frac{\|\Omega\|_{L^{\infty}}}{|x|^{n}}\chi_{\frac{3}{4}\cdot2^{k+s}\leq |x|\leq\frac{5}{4}\cdot2^{k+s}}\lesssim \frac{\|\Omega\|_{L^{\infty}}}{|x|^{n}}\chi_{
|x|\ge\frac{3}{4}\cdot2^{k}}.
\end{align} On the other hand, we compute the gradient. Again by the support of  $\varphi$, we have
\begin{align*}
\nabla H_{k,s,N(j)}(x)
&=v_{k+s}\ast\nabla\varphi_{{k+s-N(j)}}(x)\\
&=\frac{1}{2^{k+s-N(j)}}v_{k+s}\ast(\nabla\varphi)_{{k+s-N(j)}}(x)\\
&=\frac{1}{2^{k+s-N(j)}}\int_{2^{k+s}< |y|\le
2^{k+s+1}}\frac{\Omega(y)}{|y|^{n}}(\nabla\varphi)_{{k+s-N(j)}}(x-y)dy.
\end{align*}
Since $|x-y|\leq\frac{2^{k+s}}{100}, 2^{k+s}< |y|\le 2^{k+s+1}$, then $|x|\simeq |y|$ and $\frac{3}{4}\cdot2^{k+s}\leq |x|\leq\frac{5}{4}\cdot2^{k+s}.$ Thus, we get
\begin{align}\label{h2}
|\nabla H_{k,s,N(j)}(x)|
&\lesssim\|\Omega\|_{L^{\infty}}\frac{2^{N(j)}}{|x|^{n+1}}\chi_{\frac{3}{4}\cdot2^{k+s}\leq |x|\leq\frac{5}{4}\cdot2^{k+s}}.
\end{align}
Therefore, if $|h|<\frac{|x|}{2},$ we get $|x-\theta h|\simeq |x|$ and
\begin{align}\label{v7}
|K_{k,j}^\mathcal{N}(x-h)-K_{k,j}^\mathcal{N}(x)|&\leq |\nabla K_{k,j}^\mathcal{N}(x-\theta
h)||h|\\&\lesssim\sum_{s\ge 0}\|\Omega\|_{L^{\infty}}\frac{2^{N(j)}}{|x|^{n+1}}\chi_{\frac{3}{4}\cdot2^{k+s}\leq |x|\leq\frac{5}{4}\cdot2^{k+s}}(x)\cdot |h|\nonumber\\
&\lesssim\|\Omega\|_{L^{\infty}}\frac{2^{N(j)}|h|}{|x|^{n+1} } \chi_{ |x|\ge\frac{3}{4}\cdot2^k}(x)\nonumber.
\end{align}
 If $|h|<\frac{|x|}{2},$ combining \eqref{h1} and \eqref{v7}, we get for $j\ge 0$ and $k\in \Bbb Z,$
\begin{align*}
|K_{k,j}^\mathcal{N}(x-h)-K_{k,j}^\mathcal{N}(x)|
&\lesssim\frac{\omega_j(\frac{|h|}{|x|})}{|x|^n}\chi_{ |x|\ge\frac{3}{4}\cdot2^k}(x),
\end{align*} where $\omega_j(t)\le\|\Omega\|_{L^{\infty}}\min(1, 2^{N(j)}t).$

The proof of Lemma \ref{l1} is complete. \qed

We now continue the proof of \eqref{pro03}. By noting that
 supp $K_{k,j}^\mathcal{N}\subset \{x: |x|\ge\frac{3}{4}\cdot2^k\}$ and that $K_{k,j}^\mathcal{N}$ satisfies mean value zero, \eqref{g8} and \eqref{gr8}, repeating the argument of Proposition \ref{pro4} only with $K^k$ replaced by $K_{k,j}^\mathcal{N},$ we can also get
 for $1<p<\infty$ and $w\in
A_p$\begin{align}\label{I3}\bigg\|\bigg(\sum_{k\in \Bbb Z}|(\delta_0-\phi_k)\ast T_{k,j}^\mathcal{N} f|^q\bigg)^{1/q}\bigg\|_{L^p(w)}&\lesssim
\|\Omega\|_{L^{\infty}}(1+N(j))^{1+1/q}\|f\|_{L^p(w)}. \end{align}

Taking $w=1$ in \eqref{I3} we get for $1<p<\infty,$ \begin{align}\label{tq2} \bigg\|\bigg(\sum_{k\in \Bbb Z}|(\delta_0-\phi_k)\ast  T_{k,j}^\mathcal{N}
f|^q\bigg)^{1/q}\bigg\|_{L^p}&\lesssim \|\Omega\|_{L^{\infty}}(1+N(j))^{1+1/q}\|f\|_{L^p}. \end{align}
Interpolating between \eqref{phi3} and \eqref{tq2}, we get\begin{align}\label{t2} \bigg\|\bigg(\sum_{k\in \Bbb Z}|(\delta-\phi_k)\ast  T_{k,j}^\mathcal{N}
f|^q\bigg)^{1/q}\bigg\|_{L^p}&\lesssim \|\Omega\|_{L^{\infty}}2^{-\tau N(j-1)}(1+N(j))^{1+1/q}\|f\|_{L^p}. \end{align}

Similar to the proof of Theorem \ref{thm3},
based on \eqref{I3} and \eqref{t2}, applying the interpolation theorem with change of measures (Lemma \ref{sw}),  we get$$ \bigg\|\bigg(\sum_{k\in \Bbb Z}|(\delta-\phi_k)\ast  T_{k,j}^\mathcal{N} f|^q\bigg)^{1/q}\bigg\|_{L^p(w)}\lesssim
\|\Omega\|_{L^{\infty}}\{w\}_{A_p}(w)_{A_p}^{1+1/q}\|f\|_{L^p(w)},$$ which gives \eqref{pro03}.

The proof of Theorem \ref{pro01} is complete.\qed
\quad
\medskip

\section{Proof of Theorem \ref{pro02}}
Recall that we can write $$
\begin{cases}
\mathcal{S}_q( \mathcal{{T}}_\Omega f)(x)=\bigg(\sum\limits_{k\in\mathbb Z}[V_{q,k}( f)(x)]^q\bigg)^{1/q};\\[9pt]
V_{q,k}( f)(x)=\bigg(\sup\limits_
{2^k\leq t_0<\cdots<t_\lambda<2^{k+1}}{\sum\limits_{l=0}^{\lambda-1}}| T_{k,t_{l+1}}f(x)- T_{k,t_l}f(x)|^q\bigg)^{1/q},
\end{cases}
$$
where \begin{align*}T_{k,t_{l}}f(x)=\int_{t_l\le |x-y|\le 2^{k+1}}\frac{\Omega(x-y)}{|x-y|^n}f(y)\,dy.\end{align*}

{ Observe that
 \begin{align*}
 \mathcal{S}_{\infty}(\mathcal{{T}}_\Omega f)(x)\le \|\Omega\|_{L^\infty}Mf(x).
 \end{align*}
 Then by the sharp weighted boundedness of the Hardy--Littlewood maximal operator $M$ (see Hyt\"{o}nen--P\'{e}rez \cite[Corollary 1.10]{HP}, the original version
 was due to  Buckley \cite{Buc}),
 $$  \|Mf\|_{L^p(w)}\leq c_n \cdot p' \cdot  [w]_{A_p}^{1\over p} [w^{1-p'}]_{A_\infty}^{1\over p}\|f\|_{L^p(w)},\quad 1<p<\infty,$$
   we get
\begin{align}\label{s1f}
\|\mathcal{S}_{\infty}({\mathcal{T}}_\Omega f)\|_{L^p(w)}&\le c_n\cdot p'\cdot  \|\Omega\|_{L^\infty}\{w\}_{A_p}\|f\|_{L^p(w)}.
\end{align}
}
Now we claim that
\begin{align}\label{s2f}
\|\mathcal{S}_{2}(\mathcal{T}_\Omega f)\|_{L^p(w)}&\lesssim  \|\Omega\|_{L^\infty}(w)_{A_p}^{1/2}\{w\}_{A_p}\|f\|_{L^p(w)}.
\end{align}
In fact, interpolating between \eqref{s1f} and \eqref{s2f}, we get for $q\ge 2,$
\begin{align}\label{sqf}
\|\mathcal{S}_{q}({\mathcal{T}}_\Omega f)\|_{L^p(w)}&\lesssim  \|\Omega\|_{L^\infty}(w)_{A_p}^{1/q}\{w\}_{A_p}\|f\|_{L^p(w)}.
\end{align}

\medskip\medskip

Now we turn to verifying \eqref{s2f}.
Recall that  $S_{j}(f):=f\ast\varphi_{j}$ and  their difference
\begin{align}\label{eq-0s}
S_{j}(f)-S_{j+1}(f):=f\ast\psi_{j},
\end{align} where $\varphi_{j}$ and $\psi_{j}$ are
defined as  in  Section 5.
Similarly,  for any sequence of integers $\mathcal{N}=\{N(j)\}_{j=0}^{\infty}$, with
$0 = N(0) < N(1) <\cdot\cdot\cdot < N(j)\rightarrow\infty $, we also have the identity
\begin{align}\label{tn}
T_{k,t}&=T_{k,t}S_k+\sum_{j=1}^\infty T_{k,t}(S_{k-N(j)}-S_{k-N(j-1)})=:\sum_{j=0}^\infty T^\mathcal{N}_{k, t, j},
\end{align}
where $T_{k,t,0}:=T_{k,t,0}^\mathcal{N}:=T_{k,t}S_k$ and for $j\geq1$,
\begin{align}\label{tkn}
 T^\mathcal{N}_{k,t,j}&=T_{k,t}(S_{k-N(j)}-S_{k-N(j-1)})=\sum_{i=N(j-1)+1}^{N(j)}T_{k,t}(S_{k-i}-S_{k-(i-1)})=:\sum_{i=N(j-1)+1}^{N(j)}T_{k,t,i}.
\end{align}Therefore, by the Minkowski inequality, we get\begin{align}\label{s1}
\mathcal{S}_2({\mathcal{T}}_\Omega f)(x)
&=\bigg(\sum_{k\in\mathbb{Z}}|V_{2,k}(
f)(x)|^2\bigg)^\frac{1}{2}:=\bigg(\sum_{k\in\mathbb{Z}}\bigg\|\bigg\{\sum^\infty_{j=0}T_{k,t,j}^\mathcal{N}f(x)\bigg\}_{t\in[1,2]}\bigg\|_{V_2}^2\bigg)^\frac{1}{2}\\
&\leq\sum^\infty_{j=0}\bigg(\sum_{k\in\mathbb{Z}}\bigg\|\bigg\{T_{k,t,j}^\mathcal{N}f(x)\bigg\}_{t\in[1,2]}\bigg\|_{V_2}^2\bigg)^\frac{1}{2}\nonumber\\
&=:\sum^\infty_{j=0}\mathcal{S}_{2,j}^\mathcal{N}( f)(x)\nonumber.
\end{align}
Then for $1<p<\infty$ and $w\in A_p,$
\begin{align*}
\|\mathcal{S}_2({\mathcal{T}}_\Omega f)\|_{L^p(w)}\leq\sum^\infty_{j=0}\|\mathcal{S}_{2,j}^\mathcal{N}( f)\|_{L^p(w)}.
\end{align*}
\medskip
\noindent\textbf{Part 1. }We first give the $L^2$-norm of  $\mathcal{S}_{2,j}^\mathcal{N}( f)(x)$.
\begin{lemma}  \label{l1s} Let $\mathcal{S}_{2,j}^\mathcal{N}$  be defined as in \eqref{s1}. There is a positive $0<\tau<1$ such
that
\begin{align}\label{s2}
\|\mathcal{S}_{2,j}^\mathcal{N}(f)\|_{ L^2}\lesssim \|\Omega\|_{L^\infty}2^{-\tau N(j-1)}(N(j)+1))^{1/2}
\|f\|_{L^2},\quad j\geq0.
\end{align}
\end{lemma}
\emph{Proof. } For $t\in[1,2]$, we define $v_{0,t}$ as
$$v_{0,t}(x)=\frac{\Omega(x')}{|x|^n}\chi_{_{t<|x|\leq2}}(x)$$
and $v_{k,t}(x)={2^{-kn}}\nu_{0,t}(2^{-k}x)$ for $k\in\mathbb{Z}$. Then $T_{k,t}$ can be expressed as $T_{k,t}f(x)=v_{k,t}\ast f(x)$, $t\in [1,2]$.  Using \begin{align*}
\|a\|_{V_{2}}\leq \|a\|^{1/2}_{L^2}\|a'\|^{1/2}_{L^2},
\end{align*}
where $a'=\{\frac{d}{dt}a_{t}:t\in\mathbb{R}\}$  (see \cite{JSW08}), we have
\begin{align*}
\bigg(\mathcal{S}_{2,j}^\mathcal{N}( f)(x)\bigg)^2\leq\sum_{k\in\mathbb{Z}}\bigg(\int_1^2|T_{k,t,j}^\mathcal{N}f(x)|^2\frac{dt}{t}\bigg)^\frac{1}{2}
\bigg(\int_1^2|\frac{d}{dt}T_{k,t,j}^\mathcal{N}f(x)|^2\frac{dt}{t}\bigg)^\frac{1}{2}.
\end{align*}
This, along with the Cauchy--Schwartz inequality, yields
\begin{align}\label{tlm}
\|\mathcal{S}_{2,j}^\mathcal{N}(  f)\|_{L^2}^2
&\leq\bigg\|\bigg(\int_1^2\sum_{k\in\mathbb{Z}}|T_{k,t,j}^\mathcal{N}f|^2\frac{dt}{t}\bigg)^\frac{1}{2}\bigg\|_{L^2}
\bigg\|\bigg(\int_1^2\sum_{k\in\mathbb{Z}}|\frac{d}{dt}T_{k,t,j}^\mathcal{N}f|^2\frac{dt}{t}\bigg)^\frac{1}{2}\bigg\|_{L^2}.
\end{align}
By the Plancherel theorem and \eqref{tkn}, we get for $j\ge1,$
\begin{align}\label{l2}
\bigg\|\bigg(\int_1^2\sum_{k\in\mathbb{Z}}|T_{k,t,j}^\mathcal{N}f|^2\frac{dt}{t}\bigg)^\frac{1}{2}\bigg\|_{L^2}
&=\bigg(\int_1^2\sum_{k\in\mathbb{Z}}|\widehat{T_{k,t,j}^\mathcal{N}f}(\xi)|^2\frac{dt}{t}d\xi\bigg)^{1/2}\\ \nonumber
&\le\sum_{i=N(j-1)+1}^{N(j)}\bigg(\int_1^2\sum_{k\in\mathbb{Z}}|\widehat{T_{k,t,i}f}(\xi)|^2\frac{dt}{t}d\xi\bigg)^{1/2}.
\end{align}
Since $|\widehat{\psi}(2^{k-i}\xi)|\lesssim \min(|2^{k-i}\xi|,1),\,|\widehat{v}_{k,t}(\xi)|\lesssim \|\Omega\|_{L^\infty}|2^{k}\xi|^{-\beta}$ for
some $0<\beta<1$ (see \cite{DR86}),
\begin{align*}
\sum_{k\in\mathbb{Z}}|\widehat{T_{k,t,i}f}(\xi)|^2&=\sum_{k\in \mathbb{Z}}|\widehat{v}_{k,t}(\xi)|^2|\widehat{\psi}(2^{k-i}\xi)|^2|\widehat{f}(\xi)|^2\\
&\lesssim
\|\Omega\|_{L^\infty}^2\bigg(\sum_{|2^k\xi|\leq2^i}|2^{k-i}\xi|^{2}|2^{k}\xi|^{-2\beta}+\sum_{|2^k\xi|\geq2^i}|2^{k}\xi|^{-2\beta}\bigg)|\widehat{f}(\xi)|^2\\
&\lesssim\|\Omega\|_{L^\infty}^22^{-2\beta i}|\widehat{f}(\xi)|^2.
\end{align*} Then by summing the geometric series and the Plancherel theorem, we get
\begin{align}\label{l2}
\bigg\|\bigg(\int_1^2\sum_{k\in\mathbb{Z}}|T_{k,t,j}^\mathcal{N}f|^2\frac{dt}{t}\bigg)^\frac{1}{2}\bigg\|_{L^2}
&\lesssim\|\Omega\|_{L^{\infty}}2^{-\beta N(j-1)}\|f\|_{L^2}.
\end{align}
For $j=0$, $|\widehat{v}_{k,t}(\xi)|\lesssim\|\Omega\|_{L^\infty}\min (|2^k\xi|,\,|2^{k}\xi|^{-\beta})$ for some $0<\beta<1$ and
$|\widehat{\varphi}(2^{k}\xi)|\lesssim 1$,
\begin{align*}
\sum_{k\in\mathbb{Z}}|\widehat{T_{k,t,0}f}(\xi)|^2&=\sum_{k\in \mathbb{Z}}|\widehat{v}_{k,t}(\xi)|^2|\widehat{\varphi}(2^{k}\xi)|^2|\widehat{f}(\xi)|^2\\
&\lesssim\|\Omega\|_{L^\infty}^2\bigg(\sum_{|2^k\xi|\leq 1}|2^{k}\xi|^2+\sum_{|2^k\xi|\geq 1}|2^{k}\xi|^{-2\beta}\bigg)|\widehat{f}(\xi)|^2\\
&\lesssim\|\Omega\|_{L^\infty}^2|\widehat{f}(\xi)|^2.
\end{align*}
Then by  the Plancherel theorem, we get
\begin{align*}
\bigg\|\bigg(\int_1^2\sum_{k\in\mathbb{Z}}|T_{k,t,0}^\mathcal{N}f|^2\frac{dt}{t}\bigg)^\frac{1}{2}\bigg\|_{L^2}
&\lesssim\|\Omega\|_{L^{\infty}}\|f\|_{L^2}.
\end{align*}
Combining the above estimates for $j>0$ and $j=0,$ we get
\begin{align}\label{l26}
\bigg\|\bigg(\int_1^2\sum_{k\in\mathbb{Z}}|T_{k,t,j}^\mathcal{N}f|^2\frac{dt}{t}\bigg)^\frac{1}{2}\bigg\|_{L^2}
&\lesssim 2^{-\beta N(j-1)}\|\Omega\|_{L^{\infty}}\|f\|_{L^2}.
\end{align}

\quad

On the other hand,
we recall the elementary fact:  for any Schwartz function $h$,
\begin{align*}
\frac{d}{dt}[v_{j,t}\ast h(x)]&=\frac{d}{dt}\bigg[\int_{2^jt<|y|\le 2^{j+1}}\frac{\Omega(y')}{|y|^n}h(x-y)dy\bigg]\\&=\frac{d}{dt}\bigg[\int_{{\Bbb S}^{n-1}}\Omega(y')\int_{2^jt}^{2^{j+1}}\frac{1}{r}h(x-ry')drd\sigma(y')\bigg]\\
&=\frac{1}{t}\int_{{\Bbb S}^{n-1}}\Omega(y')h(x-2^jty')d\sigma(y').
\end{align*}
 Then
\begin{align}\label{df}
\bigg|\frac{d}{dt}[v_{j,t}\ast h]^{\bigwedge}(\xi)\bigg|
&=\bigg|\int_{{\Bbb S}^{n-1}}\Omega(y')e^{2\pi i2^kty'\cdot \xi}\int\frac{1}{t}h(x-2^kty')e^{2\pi i(x-2^kty')\cdot \xi}\,dxd\sigma(y')\bigg|\\ \nonumber
&=|\widehat{h}(\xi)|\ \Big|\int_{{\Bbb S}^{n-1}}\Omega(y')(e^{2\pi i2^kty'\cdot
\xi}-1)d\sigma(y')\Big|\\&\lesssim\|\Omega\|_{L^\infty}|\widehat{h}(\xi)||2^k\xi|.\nonumber
\end{align}
It is also easy to get \begin{align}\label{df1}
\bigg|\frac{d}{dt}[v_{j,t}\ast h]^{\bigwedge}(\xi)\bigg|&\lesssim\|\Omega\|_{L^\infty}|\widehat{h}(\xi)|.
\end{align}
Therefore, by the Plancherel theorem, we get for $j\ge 1$
\begin{align*}
\bigg\|\bigg(\int_1^2\sum_{k\in\mathbb{Z}}|\frac{d}{dt}T_{k,t,j}^\mathcal{N}f|^2\frac{dt}{t}\bigg)^\frac{1}{2}\bigg\|_{L^2}
&=\bigg(\int_1^2\sum_{k\in\mathbb{Z}}|\widehat{\frac{d}{dt}T_{k,t,j}^\mathcal{N}f}(\xi)|^2\frac{dt}{t}d\xi\bigg)^{1/2}\\
&\lesssim\sum_{i=N(j-1)+1}^{N(j)}\bigg(\int_1^2\sum_{k\in\mathbb{Z}}|\widehat{\frac{d}{dt}T_{k,t,i}f}(\xi)|^2\frac{dt}{t}d\xi\bigg)^{1/2}.
\end{align*}
By \eqref{df1} and $|\widehat{\psi}(2^{k-i}\xi)|\lesssim \min(|2^{k-i}\xi|, \,|2^{k-i}\xi|^{-\beta})$ for some $0<\beta<1,$\begin{align*}
\sum_{k\in\mathbb{Z}}\bigg|\widehat{\frac{d}{dt}T_{k,t,i}f}(\xi)\bigg|^2
&=\sum_{k\in\mathbb{Z}}\bigg|(\frac{d}{dt}v_{j,t}\ast\psi_{k-i}\ast
f)^{\bigwedge}(\xi)\bigg|^2\\
&\lesssim\|\Omega\|_{L^\infty}^2\sum_{k\in\mathbb{Z}}|\widehat{\psi}(2^{k-i}\xi)|^2|\widehat{f}(\xi)|^2\\&\lesssim\bigg(\sum_{|2^k\xi|\le
2^i}|2^{k-i}\xi|^2+\sum_{|2^k\xi|> 2^i}|2^{k-i}\xi|^{-2\beta}\bigg)\ \|\Omega\|_{L^\infty}^2|\widehat{f}(\xi)|^2 \\
&\lesssim \|\Omega\|_{L^\infty}^2|\widehat{f}(\xi)|^2.
\end{align*}
Then by summing the geometric series and the Plancherel theorem, we get for $j\ge 1$
\begin{align}\label{s20}
&\bigg\|\bigg(\int_1^2\sum_{k\in\mathbb{Z}}|\frac{d}{dt}T_{k,t,j}^\mathcal{N}f|^2\frac{dt}{t}\bigg)^\frac{1}{2}\bigg\|_{L^2}\lesssim
N(j)\|\Omega\|_{L^\infty}\|f\|_{L^2}.
\end{align}
For $j=0,$ by \eqref{df},\, \eqref{df1} and $|\widehat{\varphi}(2^{k}\xi)|\lesssim \min (1, \,|2^{k}\xi|^{-\beta}),$\begin{align*}
\sum_{k\in\mathbb{Z}}\bigg|\widehat{\frac{d}{dt}T_{k,t,0}^\mathcal{N}f}(\xi)\bigg|^2\big|\widehat{f}(\xi)\big|^2&=\sum_{k\in\mathbb{Z}}\bigg|\Big(\frac{d}{dt}v_{j,t}\ast\varphi_{k}\ast
f\Big)^{\bigwedge}(\xi)\bigg|^2\\
&\lesssim \|\Omega\|_{L^\infty}^2\sum_{k\in\mathbb{Z}}\min(|2^k\xi|^2,1)|\widehat{\varphi}(2^{k}\xi)|^2|\widehat{f}(\xi)|^2\\&\lesssim
\|\Omega\|_{L^\infty}^2\Big(\sum_{|2^k\xi|\le 1}|2^{k}\xi|^2+\sum_{|2^k\xi|\ge 1}|2^{k}\xi|^{-2\beta}\Big)\ |\widehat{f}(\xi)|^2 \\
&\lesssim \|\Omega\|_{L^\infty}^2 |\widehat{f}(\xi)|^2.
\end{align*}
Then by  the Plancherel theorem, we get
\begin{align}\label{s21}
&\bigg\|\bigg(\int_1^2\sum_{k\in\mathbb{Z}}|\frac{d}{dt}T_{k,t,0}^\mathcal{N}f|^2\frac{dt}{t}\bigg)^\frac{1}{2}\bigg\|_{L^2}\lesssim \|\Omega\|_{L^\infty}\|f\|_{L^2}.
\end{align}
Combining the case of $j>0$ and $j=0,$ we get \begin{align}\label{s21}
&\bigg\|\bigg(\int_1^2\sum_{k\in\mathbb{Z}}|\frac{d}{dt}T_{k,t,j}^\mathcal{N}f|^2\frac{dt}{t}\bigg)^\frac{1}{2}\bigg\|_{L^2}\lesssim(1+N(j))
\|\Omega\|_{L^\infty}\|f\|_{L^2}.
\end{align}
This along with  \eqref{tlm} and \eqref{l26}, we get for $0<\beta<1,$
\begin{align*}
&\big\|\mathcal{S}_{2,j}^\mathcal{N}( f)\big\|_{L^2}\lesssim 2^{-\frac{\beta}{2} N(j-1)}(N(j)+1)^{1/2}\|f\|_{L^2},
\end{align*} which proves Lemma \ref{l1s}.\qed

\medskip
\noindent\textbf{Part 2. }Next, we give the $L^p(w)$-norm of $\mathcal{S}_{2,j}^N( f)(x).$

For $j\ge 1,$ we denote by $S_{k,j}:=S_{k-N(j-1)}-S_{k-N(j)}$. For $j=0,$ we denote by  $S_{k,0}:=S_{k}.$  We have the following observation:
\begin{align*}\mathcal{S}_{2,j}^\mathcal{N}(
f)(x)&=\Big(\sum_{k\in\mathbb{Z}}\|\{T_{k,t}f(x)\}_{t\in[1,2]}\|_{V_2}^2\Big)^{\frac{1}{2}}\\&=\Big(\dsum_{k\in\mathbb{Z}}\dsup_{\substack
{t_1<\cdots<t_\lambda\\
[t_l,t_{l+1}]\subset[1,2]}}\dsum_{l=1}^{\lambda-1}|T_{k,t_l}
S_{k,j}f(x)-T_{k,t_{l+1}}
S_{k,j}f(x)|^2\Big)^{\frac{1}{2}}\\
&=: \Big(\dsum_{k\in\mathbb{Z}}\dsup_{\substack
{t_1<\cdots<t_\lambda\\
[t_l,t_{l+1}]\subset[1,2]}}\dsum_{l=1}^{\lambda-1}|{T}_{k,t_l,t_{l+1}}S_{k,j}
f(x)|^2\Big)^{\frac{1}{2}},\end{align*}where the operator ${{T}}_{k,t_l,t_{l+1}}$ is given by
 $${{T}}_{k,t_l,t_{l+1}}f(x):=\dint_{2^kt_l<|x-y|\le 2^{k}t_{l+1}}\frac{\Omega(x-y)}{|x-y|^n}f(y)dy=v_{k,t_l,t_{l+1}}\ast
 f(x),\,\,[t_l,t_{l+1}]\subset[1,2].$$
Denote by \begin{align}\label{gs}K_{k,l, j}:=v_{k, t_l, t_{l+1}}\ast(\varphi_{{k-N(j)}}-\varphi_{{k-N(j-1)}})\end{align} the kernel of
${T}_{k,t_l,t_{l+1}}S_{k,j }$ for $j\ge 1$ and by \begin{align}\label{gs}K_{k,l, 0}:=v_{k, t_l, t_{l+1}}\ast\varphi_{{k}}\end{align} the kernel of
${T}_{k,t_l,t_{l+1}}S_{k,0 }$ for $j=0.$
Then
\begin{align}\label{nj}
\mathcal{S}_{2,j}^\mathcal{N}(
f)(x)
&= \Bigg(\dsum_{k\in\mathbb{Z}}\dsup_{\substack
{t_1<\cdots<t_\lambda\\
[t_l,t_{l+1}]\subset[1,2]}}\dsum_{l=1}^{\lambda-1}|K_{k,l, j}\ast
f(x)|^2\Bigg)^{\frac{1}{2}},
\end{align}
In the following, we give the kernel estimates.

\quad

\begin{lemma}  \label{k1} For every $x \in \mathbb{R}^{n}\backslash\{0\}$, $j\ge 0$ and $k\in \Bbb Z,$
\begin{align}\label{v8}
\dsup_{\substack
{t_1<\cdots<t_\lambda\\
[t_l,t_{l+1}]\subset[1,2]}}\dsum_{l=1}^{\lambda-1}|K_{k,l, j}(x)|
&\lesssim\|\Omega\|_{L^{\infty}}\frac{2^k}{|x|^{n+1}}\chi_{2^{k-1}\leq |x|\leq2^{k+2}}.\end{align}If $0<|y|\le\frac{|x|}{2},$
\begin{align}\label{gr1}
\dsup_{\substack
{t_1<\cdots<t_\lambda\\
[t_l,t_{l+1}]\subset[1,2]}}\dsum_{l=1}^{\lambda-1}|K_{k,l, j}(x-y)-K_{k,l, j}(x)|
&\lesssim\frac{\omega_j(\frac{|y|}{|x|})}{|x|^n}\chi_{2^{k-1}\le|x|\le 2^{k+2}},
\end{align}where $\omega_j(t)\le\|\Omega\|_{L^\infty}\min(1, 2^{N(j)}t).$
\end{lemma}

\emph{Proof. }   Let $x \in \mathbb{R}^{n}\backslash\{0\}$. Since
supp $\varphi \subset\{x:|x| \leq \frac{1}{100}\},$ we get that
\begin{align}\label{ke}\dsup_{\substack
{t_1<\cdots<t_\lambda\\
[t_l,t_{l+1}]\subset[1,2]}}\dsum_{l=1}^{\lambda-1}|v_{k,t_l,t_{l+1}} \ast \varphi_{k-N(j)}(x)|  &=\dsup_{\substack
{t_1<\cdots<t_\lambda\\
[t_l,t_{l+1}]\subset[1,2]}}\dsum_{l=1}^{\lambda-1}\bigg|\int_{{\Bbb R}^n}\frac{\Omega(y^{\prime})}{|y|^{n}} \chi_{2^{k}t_l<|y|<2^{k}t_{l+1}} \varphi_{k-N(j)}(x-y)d
y\bigg| \\ &\lesssim\int_{{\Bbb R}^n} \frac{|\Omega(y^{\prime})|}{|y|^{n}} \chi_{2^{k}\le|y|\le2^{k+1}}
|\varphi_{k-N(j)}(x-y) |d y\nonumber\\
& \lesssim\|\Omega\|_{L^{\infty}} \frac{1}{|x|^{n}} \chi_{2^{k-1}\le|x|\le 2^{k+2}}\int_{{\Bbb R}^n}|\varphi_{k-N(j)}(x-y) | d y \nonumber\\
&\lesssim\|\Omega\|_{L^{\infty}} \frac{1}{|x|^{n}} \chi_{2^{k-1}\le|x|\le 2^{k+2}}.\nonumber \end{align}
Then we get \begin{align}\label{ke1}\dsup_{\substack
{t_1<\cdots<t_\lambda\\
[t_l,t_{l+1}]\subset[1,2]}}\dsum_{l=1}^{\lambda-1}|v_{k,t_l,t_{l+1}} \ast \varphi_{k-N(j)}(x)|
&\lesssim\|\Omega\|_{L^{\infty}} \frac{2^k}{|x|^{n+1}} \chi_{2^{k-1}\le|x|\le 2^{k+2}}. \end{align}

On the other hand, we compute the gradient. Again by taking into account the
support of $\varphi$, we obtain that
\begin{align}\label{ke2}
&\quad\dsup_{\substack
{t_1<\cdots<t_\lambda\\
[t_l,t_{l+1}]\subset[1,2]}}\dsum_{l=1}^{\lambda-1}|\nabla(v_{k,t_l,t_{l+1}}\ast \varphi_{k-N(j)})(x)|\\ &=\dsup_{\substack
{t_1<\cdots<t_\lambda\\
[t_l,t_{l+1}]\subset[1,2]}}\dsum_{l=1}^{\lambda-1}\bigg|\int_{{\Bbb R}^n} \frac{\Omega(y^{\prime})}{|y|^{n}} \chi_{2^{k}t_l<|y|<2^{k}t_{l+1}}
2^{-(k-N(j))(n+1)} \nabla \varphi\Big(\frac{x-y}{2^{k-N(j)}}\Big) d y\bigg| \nonumber\\ &\lesssim\int_{{\Bbb R}^n} \frac{|\Omega(y^{\prime})|}{|y|^{n}}\chi_{2^{k}\le|y|\le2^{k+1} } 2^{-(k-N(j))(n+1)} \bigg|\nabla
\varphi\Big(\frac{x-y}{2^{k-N(j)}}\Big)\bigg| d y\nonumber
\\ &\lesssim\|\Omega\|_{L^{\infty}} \frac{1}{|x|^{n}} \chi_{2^{k-1}\le|x|\le 2^{k+2}} \int_{{\Bbb R}^n} 2^{-(k-N(j))(n+1)}\bigg|\nabla
\varphi\Big(\frac{x-y}{2^{k-N(j)}}\Big)\bigg| d y \nonumber\\
 &\lesssim\|\Omega\|_{L^{\infty}} \frac{1}{|x|^{n}} \chi_{2^{k-1}\le|x|\le 2^{k+2}} \frac{1}{2^{k-N(j)}}\|\nabla\varphi\|_{L^1} \nonumber\\
 &\lesssim\|\Omega\|_{L^{\infty}} \frac{2^{N(j)}}{|x|^{n+1}} \chi_{2^{k-1}\le|x|\le 2^{k+2}}.\nonumber
\end{align}
Therefore, by the gradient estimate, for
$|y| \leq \frac{1}{2}|x|$ we have
\begin{align}\label{ke2}\dsup_{\substack
{t_1<\cdots<t_\lambda\\
[t_l,t_{l+1}]\subset[1,2]}}\dsum_{l=1}^{\lambda-1}|v_{k,t_l,t_{l+1}}\ast \varphi_{k-N(j)}(x-y)-v_{k,t_l,t_{l+1}}\ast \varphi_{k-N(j)}(x)| & \lesssim \frac{\|\Omega\|_{L^{\infty}}}{|x|^{n}}
2^{N(j)} \frac{|y|}{|x|}\chi_{2^{k-1}\le|x|\le 2^{k+2}}.\end{align}
Then using \eqref{ke} and \eqref{ke2},  we get
\begin{align}\label{ke3}&\dsup_{\substack
{t_1<\cdots<t_\lambda\\
[t_l,t_{l+1}]\subset[1,2]}}\dsum_{l=1}^{\lambda-1}|v_{k,t_l,t_{l+1}}\ast \varphi_{k-N(j)}(x-y)-v_{k,t_l,t_{l+1}}\ast \varphi_{k-N(j)}(x)|  \lesssim \frac{1}{|x|^{n}}
\omega_j\Big(\frac{|y|}{|x|}\Big)\chi_{2^{k-1}\le|x|\le 2^{k+2}},
\end{align}
where $\omega_j(t)\le \|\Omega\|_{L^\infty}\min(1,2^{N(j)}t)$
(for $j=0,$ the subtraction is not even needed).  From the triangle inequality,  and $N(j-1)<N(j)$ it follows that the kernel $v_{k,t_l,t_{l+1}}\ast
(\varphi_{k-N(j)}-\varphi_{k-N(j-1)})$ satisfies the same estimates \eqref{ke1} and \eqref{ke3} which proves Lemma \ref{k1}.\qed

\quad

\begin{lemma}  \label{05}   Let $\mathcal{S}_{2,j}^\mathcal{N}$ be defined as in \eqref{s1}. Then we get
\begin{align}\label{phis}\|\mathcal{S}_{2,j}^\mathcal{N}(  f)\|_{ L^{1,\infty}}\lesssim\|\Omega\|_{L^\infty}(1+N(j))^{1/2}\|f\|_{L^1}.\end{align}
\end{lemma}

\emph{{Proof.}}
We  perform
the Calder\'{o}n-Zygmund  decomposition of $f$ at height  $\alpha$,  thereby producing a disjoint  family
of dyadic cubes $\{Q\}$ with  total  measure
$$\dsum_Q|Q| \lesssim{\alpha}^{-1}\|f\|_{L^1}$$
and allowing us  to  write $f=g+h$ as in the proof of Theorem \ref{thm1}.
From Lemma \ref{l1s},
\begin{align*}
\|\mathcal{S}_{2,j}^\mathcal{N}( g)\|_{L^{1,\infty}}\lesssim \|\Omega\|_{L^\infty}(1+N(j))^{1/2}\|g\|_{L^1}.
\end{align*}
It suffices to estimate  $\mathcal{S}_{2,j}^\mathcal{N}(  h)$
away  from $\bigcup \widetilde{Q}$ where $\widetilde{Q}$ is  a   large fixed dilate  of
$Q$.  Denote
\begin{align}
\mathcal{S}_{2,j,k}^\mathcal{N} h(x)=\Big(\dsup_{\substack
{t_1<\cdots<t_\lambda\\
[t_l,t_{l+1}]\subset[1,2]}}\dsum_{l=1}^{\lambda-1}|K_{k,l, j}\ast
h(x)|^2\Big)^{\frac{1}{2}}.
\end{align}
Then from the definition of $\mathcal{S}_{2,j}^\mathcal{N}$ in \eqref{nj}
we get
\begin{align*}
\mathcal{S}_{2,j}^\mathcal{N}(h)(x)&\le \Big(\dsum_{k\in \mathbb Z}|\mathcal{S}_{2,j,k}^\mathcal{N} h(x)|^2\Big)^{1/2},
\end{align*}
which gives that
\begin{align*}
&\alpha\big|\big\{x\notin\cup \widetilde{Q}:\mathcal{S}_{2,j}^\mathcal{N} h(x)>\alpha\big\}\big|\\&\le\alpha\big|\big\{x\notin\cup \widetilde{Q}:\Big(\dsum_{k\in \mathbb Z}|\mathcal{S}_{2,j,k}^\mathcal{N} h(x)|^2\Big)^{1/2}>\alpha\big\}\big|\\&\le \alpha\big|\big\{x\notin\cup \widetilde{Q}:\dsum_{k\in \mathbb Z}|\mathcal{S}_{2,j,k}^\mathcal{N} h(x)|>(1+N(j))^{1/2}\alpha\big\}\big|+\alpha\big|\big\{x\notin\cup \widetilde{Q}:\sup_{k\in \mathbb Z}|\mathcal{S}_{2,j,k}^\mathcal{N} h(x)|>(1+N(j))^{-1/2}\alpha\big\}\big|
\\
&=: I+II.
\end{align*}

For the term $I$, by Chebyshev's inequality we have
\begin{align*}
&\alpha\big|\big\{x\notin\cup \widetilde{Q}:\dsum_{k\in \mathbb Z}|\mathcal{S}_{2,j,k}^\mathcal{N} h(x)|>(1+N(j))^{1/2}\alpha\big\}\big|\\
&\leq (1+N(j))^{-1/2} \int_{\{x\notin\cup \widetilde{Q}\}} \dsum_{k\in \mathbb Z}|\mathcal{S}_{2,j,k}^\mathcal{N} h(x)|dx\\
&\le (1+N(j))^{-1/2}\dsum_Q\int_{(\widetilde{Q})^c}\dsum_{k\in \Bbb Z}\mathcal{S}_{2,j,k}^\mathcal{N} h_Q(x)\,dx.
\end{align*}
 Denote by $y_Q$  the center of $Q.$ For $x\notin \widetilde{Q},$ we get $2|y-y_Q|\le |x-y_Q|,$ then by $|y-y_Q|\le \ell(Q)$ and using the vanishing mean value of
$h_Q$ and \eqref{gr1},
\begin{align*}
\dsum_{k\in \Bbb Z}\mathcal{S}_{2,j,k}^\mathcal{N} h_Q(x)\lesssim \dsum_{k\in \Bbb Z}\dint_{{\Bbb R}^n}|h_Q(y)|\frac{\omega_j\big(\frac{|y-y_Q|}{|x-y|}\big)}{|x-y|^n}\chi_{2^{k-1}\le
|x-y|<2^{k+2}}dy\lesssim \dint_{{\Bbb R}^n}|h_Q(y)|\frac{\omega_j\big(\frac{\ell (Q)}{|x-y|}\big)}{|x-y|^n}dy.
\end{align*}
Hence
\begin{align}\label{nq}
&\alpha\big|\big\{x\notin\cup \widetilde{Q}:\dsum_{k\in \mathbb Z}|\mathcal{S}_{2,j,k}^\mathcal{N} h(x)|>(1+N(j))^{1/2}\alpha\big\}\big|\\
&\lesssim (1+N(j))^{-1/2}\dsum_Q\dint_{{\Bbb R}^n}|h_Q(y)|\int_{|x-y|\ge 2\ell (Q)}\frac{\omega_j\big(\frac{\ell(Q)}{|x-y|}\big)}{|x-y|^n}\,dx\ dy\nonumber\\
&\lesssim (1+N(j))^{-1/2}\|\omega_j\|_{Dini}\dsum_Q\|h_Q\|_{L^1}\nonumber\\
&\lesssim  \|\Omega\|_{L^\infty}(1+N(j))^{1/2}\|f\|_{L^1},\nonumber
\end{align}
where the last inequality follows from
$$\|\omega_j\|_{Dini}\le \|\Omega\|_{L^\infty}(1+N(j)).$$
On the other hand, for the term $II$,
by \eqref{v8}, we have
\begin{align*}
\sup_{k\in \Bbb Z}|\mathcal{S}_{2,j,k}^\mathcal{N} h(x)|\le \|\Omega\|_{L^\infty}Mh(x).
\end{align*}
Then by the weak type $(1,1)$ of $M$ (see \cite{S1993}),
\begin{align}\label{nq1}
&\alpha\big|\big\{x\notin\cup \widetilde{Q}:\sup_{k\in \mathbb Z}|\mathcal{S}_{2,j,k}^\mathcal{N} h(x)|>(1+N(j))^{-1/2}\alpha\big\}\big|\lesssim \|\Omega\|_{L^\infty}(1+N(j))^{1/2}\|f\|_{L^1},
\end{align}

Combining the estimates of \eqref{nq} and \eqref{nq1}, we get
\begin{align*}
&\alpha\big|\big\{x\notin\cup \widetilde{Q}:\mathcal{S}_{2,j}^\mathcal{N} h(x)>\alpha\big\}\big|
\lesssim  \|\Omega\|_{L^\infty}(1+N(j))^{1/2}\|f\|_{L^1}.
\end{align*}

The proof of Lemma \ref{05}  is complete.\qed
\quad

\medskip
\begin{lemma}  \label{06}    For $j\ge 0,$ let $S_{2,j}^\mathcal{N}$ be defined as in \eqref{s1}. Then we get
\begin{align}\label{phi4} \mathcal{M}_{\mathcal{S}_{2,j}^\mathcal{N}}f(x)  \lesssim \mathcal{S}_{2,j}^N(
f)(x)+\|\omega_j^{2}\|_{Dini}^{1/2}Mf(x),\end{align}
where $\omega_j(t)\le\|\Omega\|_{L^\infty}\min(1, 2^{N(j)}t).$\end{lemma}

\emph{
Proof.}
 Let $Q$ be a cube, and take $x,\,\xi\in Q.$ Let $B(x)=B(x,2\sqrt{n}\ell(Q)),$ then
$3Q\subset B_x.$
By the triangular inequality, we get
\begin{align*}
& |\mathcal{S}_{2,j}^\mathcal{N}( f\chi_{{{\Bbb R}^n\setminus 3Q}})(\xi)|\\
& \le |\mathcal{S}_{2,j}^\mathcal{N}( f\chi_{{{\Bbb R}^n\setminus
B_x}})(\xi)-\mathcal{S}_{2,j}^\mathcal{N}(f{\chi_{{\Bbb R}^n\setminus B_x}})(x)|+\mathcal{S}_{2,j}^\mathcal{N}( f{\chi_{B_x\setminus
3Q}})(\xi)+\mathcal{S}_{2,j}^\mathcal{N}( f{\chi_{{\Bbb R}^n\setminus B_x}})(x)\\
&=:I+II+III.
\end{align*}

We begin by estimating term $I$.
\begin{align*}
I&\le \bigg(\sum_{k\in \Bbb Z}\bigg(\dsup_{\substack
{t_1<\cdots<t_\lambda\\
[t_l,t_{l+1}]\subset[1,2]}}\dsum_{l=1}^{\lambda-1}\int_{{\Bbb R}^n \setminus B_x}|K_{k,l, j}(\xi-y)-K_{k,l, j}(x-y)||f(y)|\,dy\bigg)^2\bigg)^{1/2}\\&\le
\Bigg(\bigg(\sum_{2^k\le \ell(Q)}+\sum_{2^k> \ell(Q)}\bigg)\bigg(\dsup_{\substack
{t_1<\cdots<t_\lambda\\
[t_l,t_{l+1}]\subset[1,2]}}\dsum_{l=1}^{\lambda-1}\int_{{\Bbb R}^n \setminus B_x}|K_{k,l, j}(\xi-y)-K_{k,l, j}(x-y)||f(y)|\,dy\bigg)^2\Bigg)^{1/2}.
\end{align*}
For $2^k\le \ell(Q),$ since $|x-\xi|\le 2|x-y|,$ we can get $|x-y|\simeq |\xi-y|.$ Therefore by \eqref{v8}
\begin{align*}
 &\sum_{2^k\le \ell(Q)}\dsup_{\substack
{t_1<\cdots<t_\lambda\\
[t_l,t_{l+1}]\subset[1,2]}}\dsum_{l=1}^{\lambda-1}\int_{{\Bbb R}^n \setminus B_x}| K_{k,l, j}(\xi-y)-K_{k,l, j}(x-y)||f(y)|\,dy\\
 &\lesssim \|\Omega\|_{L^\infty}\dsum_{2^k\le \ell(Q)}\dint_{|x-y|\ge 2\sqrt{n}\ell(Q)}\frac{2^k}{|x-y|^{n+1}}|f(y)|\,dy\\
 &\lesssim \|\Omega\|_{L^\infty}\dsum_{2^k\le \ell(Q)}2^{k-\ln\ell(Q)}Mf(x)\\
 &\lesssim \|\Omega\|_{L^\infty} Mf(x).
\end{align*}
For $2^k>\ell(Q)$, by \eqref{gr1} we get
\begin{align*}
&\bigg(\sum_{2^k> \ell(Q)}\bigg(\dsup_{\substack
{t_1<\cdots<t_\lambda\\
[t_l,t_{l+1}]\subset[1,2]}}\dsum_{l=1}^{\lambda-1}\int_{{\Bbb R}^n \setminus B_x}|K_{k,l, j}(\xi-y)-K_{k,l, j}(x-y)||f(y)|\,dy\bigg)^2\bigg)^{1/2}\\
&\lesssim  \bigg(\sum_{2^k> \ell(Q)}\bigg(\int_{|x-y|\ge 2\sqrt{n}\ell(Q)}\frac{\omega_j(|x-\xi|/|x-y|)}{|x-y|^{n}}\chi_{2^{k-1}\le |x-y|\le
2^{k+2}}|f(y)|\,dy\bigg)^2\bigg)^{1/2}\\
&\lesssim  \bigg(\sum_{2^k> \ell(Q)}\omega_j^2(\ell(Q)/2^k)\bigg)^{1/2}Mf(x)
\\&\lesssim  \|\omega_j^{2}\|_{Dini}^{1/2}Mf(x).
\end{align*}
Combining the estimates of $\sum_{2^k\le \ell(Q)}$ and $\sum_{2^k> \ell(Q)}$, we get \begin{align*}
 I\lesssim  (\|\omega^{2}\|_{Dini}^{1/2}+\|\Omega\|_{L^\infty})Mf(x).
\end{align*}

\quad

For the term $II,$  by using \eqref{v8} and noting that $|x-y|\simeq |\xi-y|$ (since $3|x-\xi|\le |x-y|$), we have
\begin{align*}
II&\lesssim \sum_{k\in \Bbb Z}\dsup_{\substack
{t_1<\cdots<t_\lambda\\
[t_l,t_{l+1}]\subset[1,2]}}\dsum_{l=1}^{\lambda-1}\int_{{\Bbb R}^n \setminus B_x}|K_{k,l, j}(\xi-y)||f{\chi_{B_x\setminus
3Q}}(y)|\,dy\\&\lesssim
\|\Omega\|_{L^\infty}\int_{3\ell(Q)\le|x-y|\le
2\sqrt{n}\ell(Q)}\frac{1}{|x-y|^n}\sum_{k\in \Bbb Z}\chi_{2^{k-1}\le |x-y|\le
2^{k+2}}|f(y)|\,dy\\&\lesssim \|\Omega\|_{L^\infty}Mf(x).
\end{align*}
\quad

We now turn to the term $III$. Note that
\begin{align*}
\mathcal{S}_{2,j}^\mathcal{N}( f{\chi_{{\Bbb R}^n\setminus B_x}})(x)
&= \Big(\dsum_{k\in\mathbb{Z}}\dsup_{\substack
{t_1<\cdots<t_\lambda\\
[t_l,t_{l+1}]\subset[1,2]}}\dsum_{l=1}^{\lambda-1}|K_{k,l, j}\ast
f{\chi_{{\Bbb R}^n\setminus B_x}}(x)|^2\Big)^{\frac{1}{2}},
\end{align*}
where for $j\ge 1,$
\begin{align*}
K_{k,l, j}\ast f(x)=v_{k,t_l,t_{l+1}}\ast (\varphi_{k-N(j)}-\varphi_{k-N(j-1)})\ast f(x)
\end{align*} and
for $j=0,$
\begin{align*}
K_{k,l, j}\ast f(x)=v_{k,t_l,t_{l+1}}\ast \varphi_{k}\ast f(x).
\end{align*}
Since  supp $v_{k,t_l,t_{l+1}}\ast \varphi_{k}\subset \{x: |x|\le 2^{k+1}\}$ and supp $v_{k,t_l,t_{l+1}}\ast\varphi_{k-N(j)}\subset \{x: |x|\le 2^{k+1}\},$
we get
\begin{align*}
III
&= \Big(\dsum_{k\in\mathbb{Z}}\dsup_{\substack
{t_1<\cdots<t_\lambda\\
[t_l,t_{l+1}]\subset[1,2]}}\dsum_{l=1}^{\lambda-1}\bigg|\int_{|x-y|>2\sqrt{n}\ell(Q)}K_{k,l, j}(x-y)f(y)\,dy\bigg|^2\bigg)^{\frac{1}{2}}\\&=
\Big(\dsum_{2\sqrt{n}\ell(Q)<2^{k+1}}\dsup_{\substack
{t_1<\cdots<t_\lambda\\
[t_l,t_{l+1}]\subset[1,2]}}\dsum_{l=1}^{\lambda-1}\bigg|\int_{|x-y|>2\sqrt{n}\ell(Q)}K_{k,l, j}(x-y)f(y)\,dy\bigg|^2\bigg)^{\frac{1}{2}}\\&\le
\bigg(\dsum_{2\sqrt{n}\ell(Q)<2^{k+1}}\dsup_{\substack
{t_1<\cdots<t_\lambda\\
[t_l,t_{l+1}]\subset[1,2]}}\dsum_{l=1}^{\lambda-1}\bigg|\int_{{\Bbb R}^n} K_{k,l, j}(x-y)f(y)\,dy
\bigg|^2\bigg)^{\frac{1}{2}}\\&\quad
+\bigg(\dsum_{2\sqrt{n}\ell(Q)<2^{k+1}}\dsup_{\substack
{t_1<\cdots<t_\lambda\\
[t_l,t_{l+1}]\subset[1,2]}}\dsum_{l=1}^{\lambda-1}\bigg|\int_{|x-y|\le 2\sqrt{n}\ell(Q)}K_{k,l, j}(x-y)f(y)\,dy\bigg|^2\bigg)^{\frac{1}{2}}
\\&\lesssim  \mathcal{S}_{2,j}^\mathcal{N}( f)(x)+\Bigg(\dsum_{2\sqrt{n}\ell(Q)<2^{k+1}}\bigg(\int_{|x-y|\le 2\sqrt{n}\ell(Q)}  \dsup_{\substack
{t_1<\cdots<t_\lambda\\
[t_l,t_{l+1}]\subset[1,2]}}\dsum_{l=1}^{\lambda-1}|K_{k,l, j}(x-y)||f(y)|\,dy\bigg)^2\Bigg)^{1/2}.
\end{align*}
Moreover, from \eqref{v8} we obtain that
\begin{align*}
&
\Bigg(\dsum_{2\sqrt{n}\ell(Q)<2^{k+1}}\bigg(\int_{|x-y|\le 2\sqrt{n}\ell(Q)}  \dsup_{\substack
{t_1<\cdots<t_\lambda\\
[t_l,t_{l+1}]\subset[1,2]}}\dsum_{l=1}^{\lambda-1}|K_{k,l, j}(x-y)||f(y)|\,dy\bigg)^2\Bigg)^{1/2}\\
 &\lesssim  \|\Omega\|_{L^\infty}\sum_{2\sqrt{n}\ell(Q)\le 2^{k+1} }2^{-kn}\int_{|x-y|\le 2\sqrt{n}\ell(Q) } |f(y)|\,dy\\&\lesssim
 \|\Omega\|_{L^\infty}Mf(x),
\end{align*}
which shows that
\begin{align*}
III&\lesssim \mathcal{S}_{2,j}^\mathcal{N}( f)(x)+\|\Omega\|_{L^\infty}Mf(x).
\end{align*}
Combining the estimates of $I,\,II$ and $III$, we get \begin{align*}
 \mathcal{S}_{2,j}^\mathcal{N}( f{\chi_{{\Bbb R}^n\setminus 3Q}})(\xi)\lesssim\mathcal{S}_{2,j}^\mathcal{N}(
 f)(x)+(\|\Omega\|_{L^\infty}+\|\omega^2\|_{Dini}^{1/2})Mf(x).
\end{align*}  which leads to  \begin{align*}
\mathcal{M}_{\mathcal{S}_{2,j}^\mathcal{N}}f(x)\lesssim \mathcal{S}_{2,j}^\mathcal{N}f(x)+(\|\Omega\|_{L^\infty}+\|\omega^2\|_{Dini}^{1/2})Mf(x).
\end{align*}

The proof of Lemma \ref{06} is complete. \qed

\quad

Then by Lemma \ref{05}, Lemma \ref{06} and the weak type $(1,1) $ of $M$ (see \cite{S1993}) we get   \begin{align}\label{vqm}
\|\mathcal{M}_{\mathcal{S}_{2,j}^\mathcal{N}}f\|_{ L^{1,\infty}} &\lesssim \|\Omega\|_{L^{\infty}}(1+N(j))^{1/2}\|f\|_{L^1}
\end{align}
since $$\|\omega^{2}\|_{Dini}^{1/2}\lesssim \|\Omega\|_{L^{\infty}}(1+N(j))^{1/2}.$$

\quad

 Since $\mathcal{S}_{2,j}^\mathcal{N}$ satisfies \eqref{phis} and \eqref{vqm}, therefore by applying Lemma \ref{lem 4} to $U=\mathcal{S}_{2,j}^\mathcal{N}$,
 we get that for $1<p<\infty$ and $w\in A_p$,
 \begin{align}\label{vs2}
  \|\mathcal{S}_{2,j}^\mathcal{N}(   f)\|_{L^p(w)}&\lesssim
 \|\Omega\|_{L^{\infty}}(1+N(j))^{1/2}\{w\}_{A_p}\|f\|_{L^p(w)}.
 \end{align}
Taking $w=1$ in the above inequality gives
 \begin{align}\label{ts2}
 \|\mathcal{S}_{2,j}^\mathcal{N}(   f)\|_{L^p}&\lesssim
\|\Omega\|_{L^{\infty}}(1+N(j))^{1/2}\|f\|_{L^p}.
\end{align}
Interpolating between \eqref{s2} and \eqref{ts2} shows that there exists some $\theta\in (0,1)$\begin{align}\label{st2} \|\mathcal{S}_{2,j}^\mathcal{N}( f)\|_{L^p}&\lesssim
\|\Omega\|_{L^{\infty}}2^{-\theta N(j-1)}(1+N(j))^{1/2}\|f\|_{L^p}. \end{align}
Combining \eqref{vs2} and \eqref{st2}, and by using  the interpolation theorem with change of measures-Lemma \ref{sw},  
we obtain that
$$
\|\mathcal{S}_{2}(\mathcal{T}_\Omega f)\|_{L^p(w)}\lesssim \|\Omega\|_{L^\infty}\{w\}_{A_p}(w)_{A_p}^{1/2}\|f\|_{L^p(w)},
$$which gives the proof of \eqref{s2f}.

The proof of Theorem 1.4 is complete.
\qed

 \end{document}